\documentclass[%
preprint,
 amsmath,amssymb,
 aps,
]{revtex4-1}

\usepackage{graphicx}
\usepackage{dcolumn}% Align table columns on decimal point
\usepackage{bm}% bold math
\usepackage{enumerate}

\usepackage{epstopdf}
\usepackage{hyperref}

\newcommand{\expected}[1]{\mathbb{E}\left[#1\right]}

\newcommand{\sgn}[1]{\operatorname{sgn}\left(#1\right)}
\renewcommand{\bar}[1]{\overline{#1}}
\newcommand{\dd}{\partial}

\begin{document}

\title{Reduced $\alpha$-stable dynamics for multiple time scale systems forced with correlated additive and multiplicative Gaussian white noise}
\author{William F. Thompson}\affiliation{Dept. of Mathematics, University of British Columbia} \email{william.thompson@alumni.ubc.ca}
\author{Rachel A. Kuske}\affiliation{School of Mathematics, Georgia Institute of Technology}\email{rachel@math.gatech.edu}
\author{Adam H. Monahan}\affiliation{School of Earth and Ocean Science, University of Victoria}\email{monahana@uvic.ca}
\date{\today}

\begin{abstract}
Stochastic averaging problems with Gaussian forcing have been studied thoroughly for many years, but far less attention has been paid to problems where the stochastic forcing has infinite variance, such as an $\alpha$-stable noise forcing. It has been shown that simple linear processes driven by correlated additive and multiplicative (CAM) Gaussian noise, which emerge in the context of atmosphere and ocean dynamics, have infinite variance in certain parameter regimes.

In this paper, we study a stochastic averaging problem where a linear CAM noise process in a particular parameter regime is used to drive a comparatively slow process. It is shown that the slow process exhibits properties consistent with being forced by a white $\alpha$-stable noise in the case of large time-scale separation. We identify the conditions required for the fast linear CAM process to have such an influence in driving a slower process, and then derive an (effectively) equivalent fast, infinite-variance process for which an existing stochastic averaging approximation is readily applied. These results are illustrated using a set of representative numerical results.
\end{abstract}

\maketitle

\section{Introduction}

%% LAYOUT OF INTRO

% What are we talking about when we talk about \alpha-stable dynamics?
Stochastic differential equation (SDE) modelling techniques are used broadly in modern quantitative research in disciplines including, but not limited to, physics \cite{Gardiner1985}, finance \cite{Cont2004}, and biology \cite{Bressloff2014}. They allow researchers to incorporate elements into dynamical models that are impossible or impractical to model explicitly due to their unpredictable nature and/or complexity, by  representating their effects with stochastic processes. The most commonly used stochastic driving process is Gaussian white noise. The ubiquity of Gaussian white noise follows from the Central Limit Theorem (CLT), which states that the sum of a sufficiently large set of independent, identically-distributed random variables with finite variance converges in distribution to a Gaussian random variable. However, the assumptions necessary for the CLT are not always satisfied. For example, if the random variables  have a density with power law tails such that the variance is not finite, then we must consider the Generalized Central Limit Theorem (GCLT) \cite{Feller1966II}. The GCLT states that a sum of $n$ independent, identically-distributed random variables of this type, not necessarily having finite variance, converges in distribution to an $\alpha$-stable random variable as $n \to \infty$. 
The distribution $\mathcal{S}_{\alpha}(\beta,\sigma)$ of an $\alpha$-stable random variable depends on three parameters: the stability index $\alpha \in (0,2]$, the skewness parameter $\beta \in [-1, 1]$, and the scale parameter $\sigma \in (0,\infty)$. Such distributions do not have in general a closed-form expression for their probability density functions (PDFs), but their characteristic functions have the form
\begin{equation}
\psi(k) = \exp\left[-\sigma^{\alpha}|k|^{\alpha}\Xi(k;\alpha,\beta) \right], \label{eqn:stable_characteristic_function}
\end{equation}
where 
\begin{align}
\Xi(k; \alpha,\beta) &= 1 - i\beta\sgn{k}\varphi(k), \label{eqn:Xi} \\
\varphi(k) &= \begin{cases}\tan(\pi\alpha/2) &\mbox{if $\alpha \ne 1$} \\ -\frac{2}{\pi}\log(|k|) &\mbox{if $\alpha = 1.$} \end{cases}
\end{align}
The case where $\alpha = 2$ is the only case of an $\alpha$-stable random variable with finite variance and corresponds to a Gaussian random variable with mean 0 and variance $2\sigma^{2}$.

\begin{figure}
\centering
\includegraphics[width=0.49\textwidth]{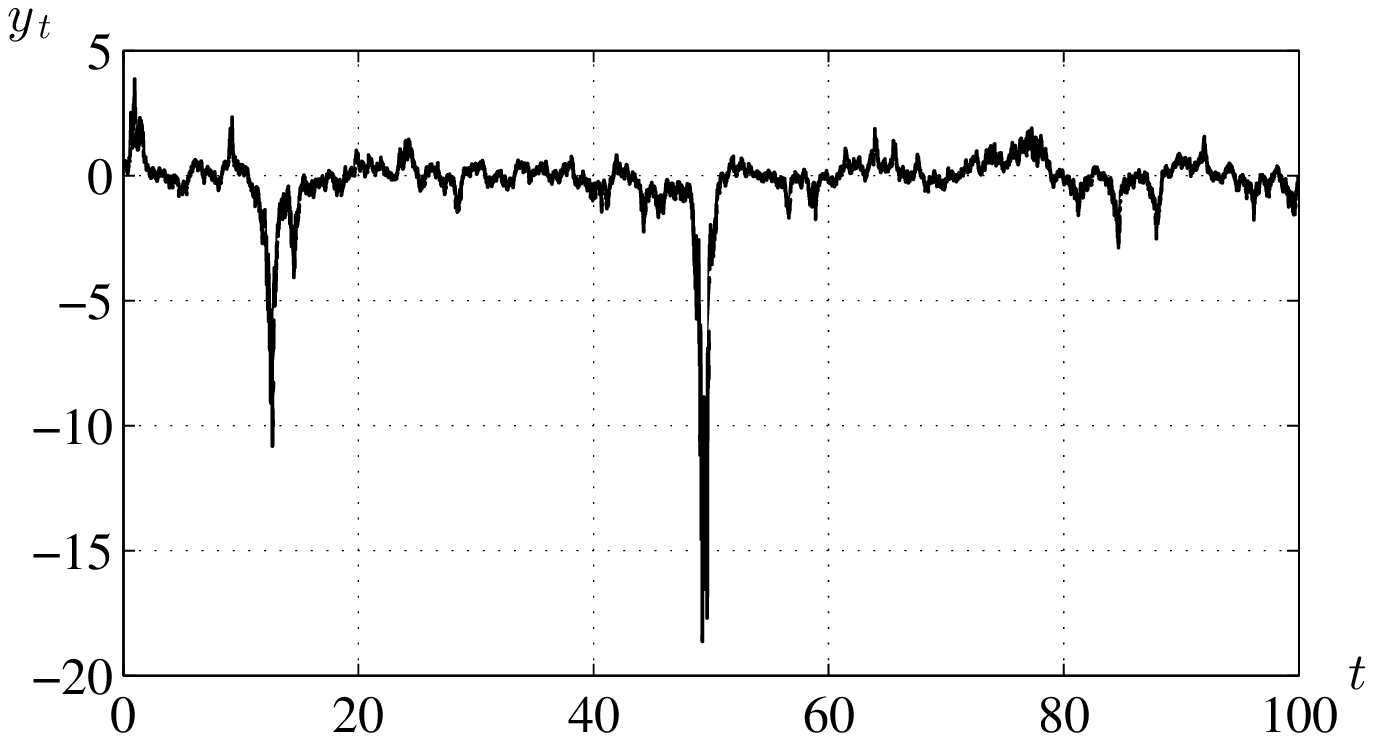}
\includegraphics[width=0.49\textwidth]{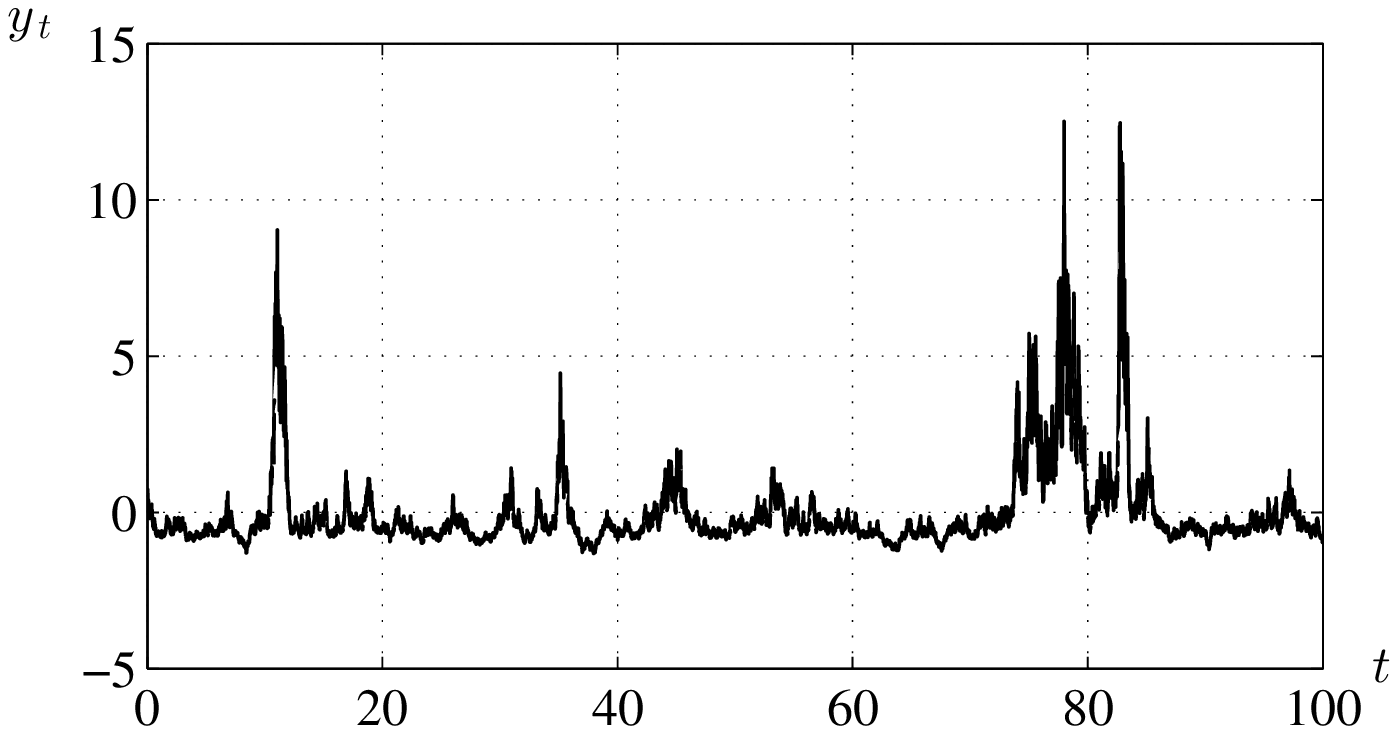}
\caption{Sample time series of $y_{t}$. {Left}: $(L,E,g,b) = (-1,1.0541,-0.1,0.5)$ corresponding to ($\alpha^{*}, \beta^{*}) = (1.8,-0.25)$. {Right}: $(L,E,g,b) = (-1,1.118,1,0.3)$ corresponding to ($\alpha^{*}, \beta^{*}) = (1.6,0.996)$.}
\label{fig:time_ser}
\end{figure}

In this paper, we consider stochastic dynamical systems with multiple time scales of the form 
\begin{align}
dx_{t} &= f_{1}(x_{t})\,dt + \epsilon^{-\rho}f_{2}(x_{t})y_{t/\epsilon}\,dt,  \label{eqn:slowSDE} \\
dy_{t} &= \left(L + {{E}^{2}}/{2}\right)y_{t}\,dt + (Ey_{t} + g)\,dW_{1,t} + b\,dW_{2,t}, \label{eqn:CAMsde}
\end{align}
where $t \ge 0$, $0 < \epsilon \ll 1$, $\rho,E,g,b$,  and $L < 0$ are real non-zero constants, 
%({\bf can $\rho =0$?}) AM: YES
 $f_{1}, f_{2}$ are functions with $f_{2}(\tilde{x}) \ne 0$ for any $\tilde{x}$ in the domain of $x_{t}$, and $dW_{1,t},\,dW_{2,t}$ are independent Gaussian white noise processes with $\expected{dW_{p,t}} = 0$, $\expected{dW_{p,s}\,dW_{q,t}} = \mathbf{1}_{p = q}\delta(t - s)dt$ for $p,q = 1,2$.  The process $y_t$ is referred to as
a CAM noise process \cite{Sardeshmukh2009, Penland2012}, described by a SDE of standard form with linear drift term and correlated additive and multiplicative (CAM) Gaussian white noise forcing.  Such processes emerge naturally in highly-truncated projections of fluid mechanical systems \cite{Sardeshmukh2009}.  The variable $y_{t/\epsilon}$ is a fast version of the process $y_{t}$, specifically satisfying
\begin{equation}
dy_{t} = \frac{1}{\epsilon} \left(L + {{E}^{2}}/{2}\right)y_{t}\,dt + \frac{1}{\epsilon^{1/2}}\left[\, (Ey_{t} + g)\,dW_{1,t} + b\,dW_{2,t}\, \right]\, .
\end{equation}
Thus the system (\ref{eqn:slowSDE})-(\ref{eqn:CAMsde}) represents a slow process $x$ coupled to a fast process $y_{t/\epsilon}$ reduced from a larger nonlinear system.  In this study, we restrict our attention to  $g\neq 0$, noting that for $g=0$ the additive and multiplicative noise is not correlated. As shown in \cite{Sardeshmukh2009}, the stationary distribution for $y_{t}$ is non-Gaussian, and has infinite variance if $-L \leq E^{2}$. Example time series for $y_{t}$ are plotted in Figure \ref{fig:time_ser} to illustrate the variability that this process can exhibit in this parameter regime.  One feature to notice is that the distribution of fluctuations in the realization are increasingly asymmetric around zero $y_t$ for larger values of $g>0$.  This skewness is discussed further in Section \ref{sec:CAMnoise}.

Our goal is to derive a one-dimensional stochastic dynamical system for a variable $X_{t}$ that weakly approximates $x_{t}$ for parameter ranges where $y_{t}$ does not possess second-order moments. Such a stochastic process $X_{t}$ is referred to as a {stochastic averaging approximation} (or {stochastic homogenization approximation}) for $x_{t}$. Stochastic averaging techniques offer systematic ways of reducing a dynamical system with multiple time scales to one with reduced dimension and fewer time scales whose properties weakly approximate those of the slow variable(s) of the full system. The prototypical stochastic averaging problem is to determine a SDE model governing the evolution of random variable $X_{t}$, such that $X_{t}$ weakly approximates $x_{t}$ where $dx_{t} = f(x_{t},z_{t/\epsilon})\,dt$, $0 < \epsilon \ll 1$, and the process $z_{t}$ evolves on a faster time scale relative to $x_{t}$.  The case where $z_{t}$ is driven by Gaussian white noise and has finite moments has been well-studied, with several established results for stochastic averaging approximations \cite{Khasminskii1966a, Khasminskii1966b, Papanicolaou1974, Borodin1977, Freidlin1984, Givon2004, Pavliotis2007}.  Some studies consider a fast chaotic driving process rather than a stochastic one \cite{Kantz2004, Mitchell2012}.  The case where the fast variable is driven by an $\alpha$-stable noise process has received less attention, with some stochastic averaging approximations obtained in \cite{Srokowski2011,Thompson2014}.  
 
Stochastic averaging techniques have particular importance to climate modelling problems characterized by high dimensionality, chaotic dynamics, and multiple time scales \cite{Hasselmann1976, Majda2001, Arnold2003, Mona11, Mitchell2012}. The multiple time scale nature of climate variability is two-fold. Climate variability involves interactions between different Earth system components with broadly separated dominant timescales (e.g. atmosphere, ocean, cryosphere, land surface). Furthermore, each of these components are high-dimensional systems with a broad range of timescales of internal variability \cite{Saltzman2002}. Linear CAM noise processes of the form (\ref{eqn:CAMsde}) have been used %in a number of studies 
 to understand certain various features of the probability distributions of observed climate variables \cite{Sura2004, Penland2008, Sura2008, Sardeshmukh2009, Sura2011}.  Since infinite variance processes can arise in simple systems like (\ref{eqn:CAMsde}) with additive and multiplicative noise,  the GCLT suggests that CAM noise dynamics could be a potential source of the appearance of $\alpha$-stable noise forcing in the long-term climate record \cite{Ditlevsen1999}.   However, the results in \cite{Penland2012} demonstrate that linear CAM noise processes are not equivalent to $\alpha$-stable noise processes, so that the recent results of \cite{Thompson2014} can not be applied directly to (\ref{eqn:slowSDE})-(\ref{eqn:CAMsde}).  Rather,  the GCLT indicates that similarities
 in the behavior of the time integrals of the CAM noise and $\alpha$-stable processes can be demonstrated over sufficiently long time scales.

%From the GCLT, a sum of $n$ independent CAM processes will converge to an alpha stable process as $n$ increases. << DO I WANT TO ADD THIS SENTENCE? >>

The main result of this paper is to identify a system composed of $x_t$ in 
(\ref{eqn:slowSDE}) coupled with a Ornstein-Uhlenbeck-L\'{e}vy process (OULP) $z_t$ that, under a stochastic averaging approximation, gives a slow proess $X(t)$ that weakly approximates the slow dynamics of $x_t$ in (\ref{eqn:slowSDE})-(\ref{eqn:CAMsde}). Analogous to an Ornstein-Uhlenbeck (OU) process forced by Gaussian noise, the OULP $z_{t/\epsilon}$ is a one-dimensional process with a linear drift driven by an $\alpha$-stable white noise forcing $dL_{t}^{(\alpha,\beta)} \sim \mathcal{S}_{\alpha}(\beta,dt^{1/\alpha})$.
 The key to this main result is a comparison of the integral of $y_t$ in (\ref{eqn:CAMsde}) with the integral of $z_t$.
Once such a process $z_t$ is identified, the (N+) approximation for systems driven with additive $\alpha$-stable noise given in \cite{Thompson2014} can be used to derive a stochastic averaging approximation for the slow dynamics. 

The remainder of the paper is organized as follows. In Section \ref{sec:CAMnoise}, we give some properties of the linear CAM process $y_{t}$ and discuss its relationship with $\alpha$-stable distributions through the GCLT. In Section \ref{sec:reductionCAMtostable}, we identify an 
 OULP $z_{t/\epsilon}$  whose integral weakly approximates that of  $y_{t/\epsilon}$ in (\ref{eqn:CAMsde}). Replacing $y_{t/\epsilon}$ with $z_{t/\epsilon}$, in Section \ref{sec:stochastic_averaging_approx_final} we derive the stochastic averaging approximation for the slow variable $x_t$ in (\ref{eqn:slowSDE}) by applying the (N+) stochastic averaging approximation as in
\cite{Thompson2014}. 
  As a minimum requirement, the characteristic time scale of $y_{t/\epsilon}$ must be at least an order of magnitude faster than $x_{t}$ for $X_{t}$ to provide an accurate weak approximation to $x_{t}$. However the exact time scale separation needed depends on the parameters of the CAM noise process, as the rate of convergence to the limiting $\alpha$-stable distribution depends on the tail behaviour of the random variables being summed in the integral of $y_t$.
 In Section \ref{sec:sample_systems}, we apply our approximation to specific examples, one linear system and two nonlinear systems of the form (\ref{eqn:slowSDE})-(\ref{eqn:CAMsde}), illustrating the performance of the approximation. 
 
\section{The CAM noise process and its relationship to $\alpha$-stable distributions}
\label{sec:CAMnoise}
In this section, we state some basic properties of the CAM noise process $y_t$ (many of which are derived in \cite{Sardeshmukh2009, Penland2012}), identify the parameter domains of interest, and discuss the relationship between the CAM noise distribution and the $\alpha$-stable distribution for these parameters.
\subsection{Properties of the linear CAM noise process}
\label{subsec:properties_of_CAM}
 First we note that $y_{t}$ can display a range of different behaviours, despite its relatively simple appearance. It reduces to an Ornstein-Uhlenbeck process (OUP) as $E \to 0$ and to a geometric Brownian motion when both $E/g,\,E/b \to \infty$. We focus on parameter ranges for which
$y_t$ has a stationary probability density function without a finite variance. The relevant parameter ranges can be determined by considering the probability density function (PDF) $p(y,t)$ for (\ref{eqn:CAMsde}), which satisfies the FKE
\begin{align}
\frac{\dd p}{\dd t} = \mathcal{A}p, \quad \mathcal{A}p \nonumber &= -\left(L + \frac{E^{2}}{2}\right)\frac{\dd}{\dd y}(yp) \\ &+ \frac{1}{2}\frac{\dd^{2}}{\dd y^{2}}\left([(Ey + g)^{2} + b^{2}] p\right). \label{eqn:CAMfke}
\end{align}
The stationary PDF, $p_s(y)$ is obtained by solving the time-independent FKE, $\mathcal{A}p_s = 0$. This function is given in \cite{Penland2012},
\begin{align}
p_s(y) = &\frac{1}{\mathcal{N}}\left[(Ey + g)^{2} + b^{2} \right]^{-(\nu + 1)} \nonumber \\ &\times \exp\left(\frac{2g\nu}{b}\arctan\left(\frac{Ey + g}{b}\right) \right),\label{eqn:CAMpdf}
\end{align}
where
\begin{equation}
\nu = -\left(\frac{L}{E^{2}} + \frac{1}{2}\right) \label{nueq}
\end{equation}
and $\mathcal{N}$ is the normalization constant, shown in \cite{Penland2012} to be given by
\begin{align}
\mathcal{N} &= \frac{1}{b^{2\nu + 1}E}\int_{-\pi/2}^{\pi/2}\frac{\exp\left( 2g\nu\xi/b\right)}{(1 + \tan^{2}(\xi))^{\nu}}\,d\xi \\ &= \frac{2\pi(2b)^{-(2\nu + 1)}\Gamma(2\nu + 1)}{E\Gamma(\nu + 1 - i\frac{2g\nu}{b})\Gamma(\nu + 1 + i\frac{2g\nu}{b})}
\label{eqn:CAM_normalization_cons}
\end{align}
where $\Gamma$ is the complex Gamma function. If $\nu \leq 0$, the first moment of the distribution (\ref{eqn:CAMpdf}) does not exist. We consider $\nu > 0$ for the remainder of this paper, and hence the stationary mean of $y_{t}$ is equal to 0. The distribution is heavy-tailed \cite{Penland2012}, since for large $|y|$,
 $p_s$ (\ref{eqn:CAMpdf}) decays according to a power law,  
\begin{equation}
p_s(y) \sim \frac{h(\sgn{y})}{|y|^{2(\nu + 1)}} \quad \mbox{as} \quad
|y| \rightarrow \infty, \label{eqn:tail_behaviour}
\end{equation}
where $h(s) = \frac{\exp\left(\frac{\pi g\nu}{b}s \right)}{\mathcal{N}E^{2(\nu + 1)}}$. Eq. (\ref{eqn:tail_behaviour}) indicates that $y$ with stationary density $p_s$ in  (\ref{eqn:CAMpdf}) does not have a finite variance $0<\nu \leq \frac{1}{2}$. By the GCLT, this case is related to non-Gaussian $\alpha$-stable processes and so we restrict ourselves to corresponding values of
$L$ and $E$ in (\ref{nueq}), that is, $\frac{E^{2}}{2} < (-L) < E^{2}$. 

\begin{figure}[h]
\centering
\includegraphics[width=0.49\textwidth]{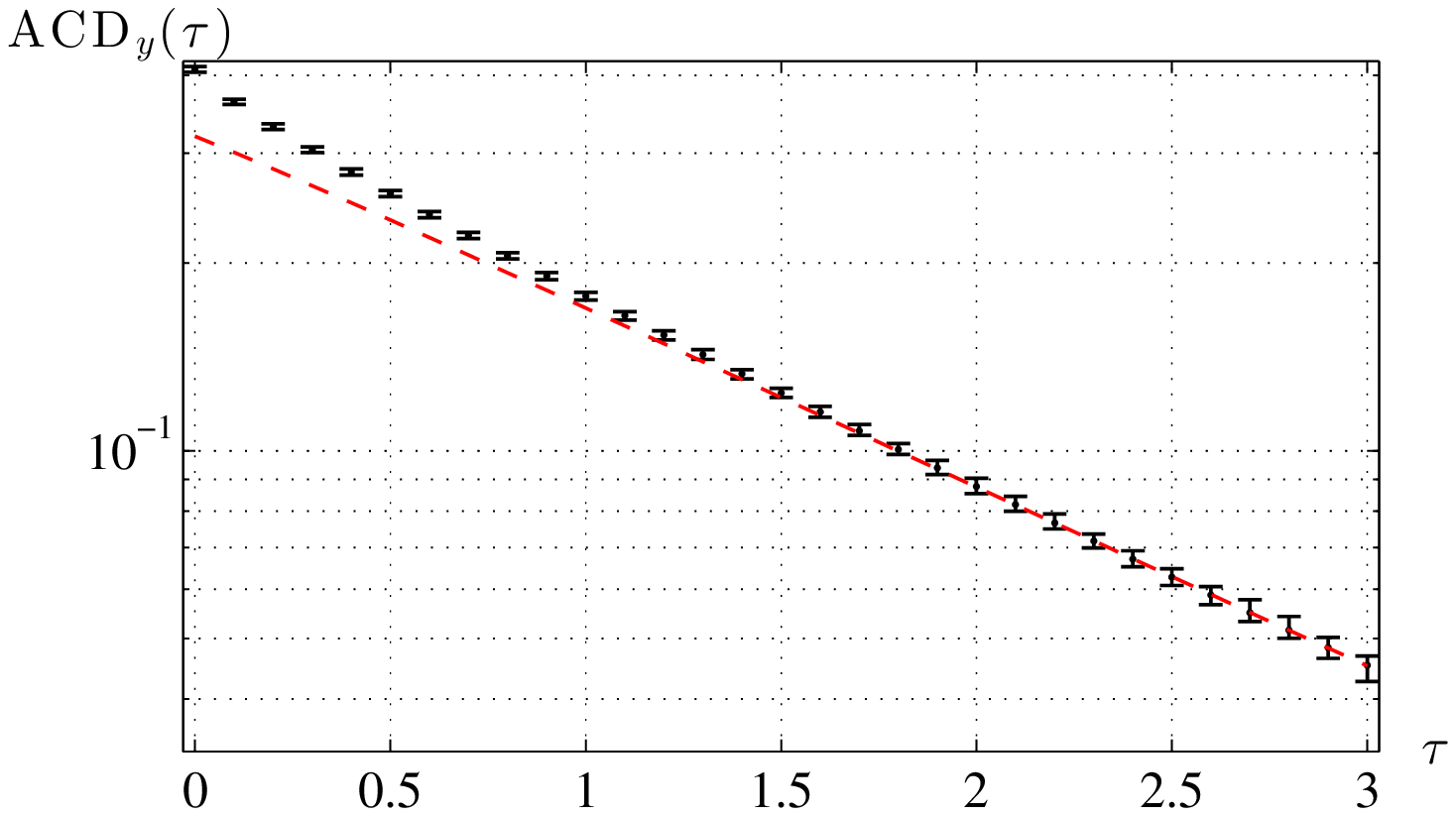}
\includegraphics[width=0.49\textwidth]{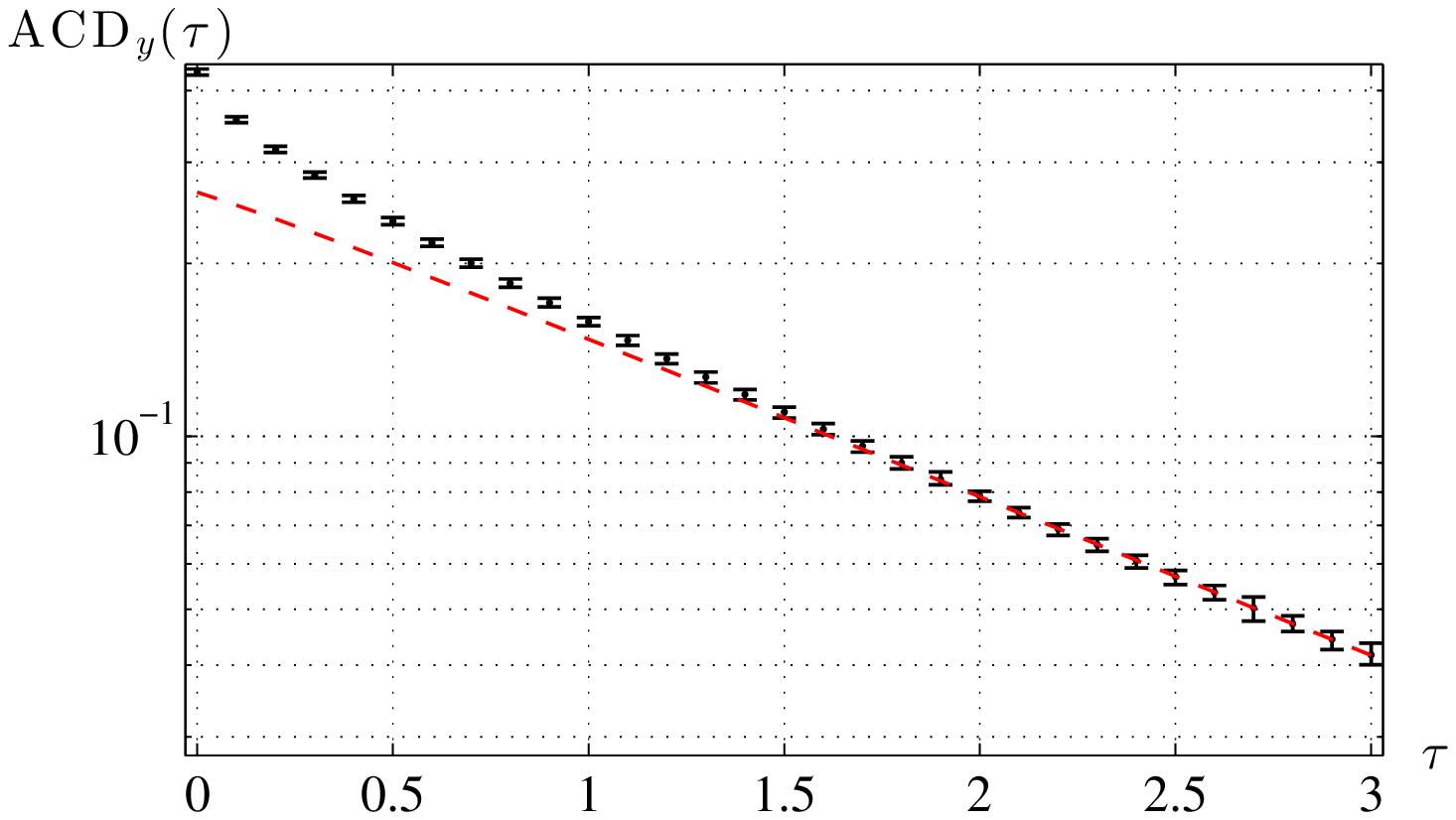}
\caption{Logarithmic plots of numerical estimates of the real part of $\operatorname{ACD}_{y}(T)$ for the CAM processes shown
in Figure \ref{fig:time_ser}:  {Left}: $(L,E,g,b) = (-1,1.0541,-0.1,0.5)$  corresponding to ($\alpha^{*}, \beta^{*}) = (1.8,-0.25)$; {Right}: $(L,E,g,b) = (-1,1.118,1,0.3)$ corresponding to ($\alpha^{*}, \beta^{*}) = (1.6,0.996)$.  The ACD is estimated from 100 realizations from $t = 0$ to $t = 40000$.. The imaginary part of $\operatorname{ACD}_{y}(T)$ is  several orders of magnitude smaller and is not shown. Error bars indicate the 25th and 75th percentiles of the estimates of $\operatorname{ACD}_{y}(T)$ across the different realizations. The red dashed line indicates the choice of $\theta$ in $\operatorname{ACD}_{z}(T)$ (\ref{eqn:autocodiff_linear}) that approximates the slope of $\log\left(\operatorname{ACD}_{y}(T)\right)$ for larger values of $\tau$.}
\label{fig:sample_acds}
\end{figure}

To capture the serial dependence of $y_{t}$ (or, less formally, the memory of $y_{t}$)  we use the autocodifference function (ACD function),   the quantity for stochastic processes without a finite variance analogous to the autocovariance for processes with
finite variance. The ACD for  the stationary process $y_t$
% with $E \neq 0$
 is defined in terms of characteristic function of the process at different times \cite{Taqqu1994, Wylomanska2015}: 
\begin{equation}
\operatorname{ACD}_{y}(\tau) = \log\left[\frac{\expected{\exp(i(y_{t+\tau} - y_{t}))}}{\expected{\exp(iy_{t+\tau})}\expected{\exp(-iy_{t})}}\right]. \label{eqn:autocodiff}
\end{equation}
Figure \ref{fig:sample_acds} presents sample estimates of (\ref{eqn:autocodiff}) for CAM processes corresponding to two sets of parameter values.  For comparison, also shown are ACDs for the OULPs $z(t)$
\begin{eqnarray}
dz_{t} &=& -\theta z_{t}\,dt + \sigma_z\,dL_{t}^{(\alpha,\beta)}, \quad \theta > 0, \\ 
dL_{t}^{(\alpha,\beta)} &\sim& \mathcal{S}_{\alpha}(\beta,dt^{1/\alpha}). \label{eqn:zt_dynamics} \\
\operatorname{ACD}_{z}(\tau) &=& \frac{\sigma_z^{\alpha}}{\alpha\theta}\bigg\{1 + \exp(-\alpha\theta \tau) - |1 - \exp(-\theta \tau)|^{\alpha} -  \nonumber \\
&&\left. i\beta\tan\left(\frac{\pi\alpha}{2}\right)\left[ (1 - \exp(-\alpha\theta \tau)) - |1 - \exp(-\theta \tau)|^{\alpha}\right]\right\}. \label{eqn:autocodiff_linear}
\end{eqnarray}
The parameter values of these OULPs were selected so that the asymptotic decay rates of their ACD functions match those of the corresponding CAM processes.
We observe that  $\log[ {\rm Re}(\operatorname{ACD}_y(\tau))]$ of the CAM noise process is evidently a nonlinear function of the lag $\tau$, in contrast to  $\log[{\rm Re}(\operatorname{ACD}_z)]$ for the OULP which is close to being linear in $\tau$. Specifically, for shorter lag times $\tau$, $\operatorname{ACD}_y(\tau)$ decays more rapidly than the $\operatorname{ACD}_z$ for the  OULP with the same long-time dependence structure, while for larger values of $\tau$, $\log[{\rm Re}(\operatorname{ACD_y)]}$ asymptotes to approximately linear behavior. This serial dependence plays an important role in considering the asymptotic behavior of integrals of $y_t$ in 
Section \ref{sec:reductionCAMtostable} below.

\subsection{Generalized central limit theorem applied to independent CAM random variables} 
\label{sec:CAMandStable}
The GCLT states that a scaled and shifted sum of independent, identically distributed (iid) random variables $R_{j}$ with density $u_{R}(r) \propto |r|^{-(\alpha + 1)}$ as $|r| \rightarrow \infty$ with $\alpha \in (0,2)$,
 converges in distribution to an $\alpha$-stable random variable (rather than a Gaussian random variable) \cite{Feller1966II},  denoted
\begin{equation}
\frac{1}{n^{1/\alpha}}\sum_{j = 1}^{n}\left( R_{j} - \bar{R}_{n} \right) \quad \underset{D}{\rightarrow} \quad \mathcal{S}_{\alpha^*}(\beta^*,\sigma^*) \quad {\mbox{ for } }\alpha = \alpha^*\, .\label{eqn:CAM-to-stable}
\end{equation}
We illustrate how the parameters $\alpha^*$, $\beta^*$, and $\sigma^*$ are determined in the context where the independent random variables $\{R_{j}\}_{j = 1}^{n}$ are drawn from the  distribution $p_s$ given by (\ref{eqn:CAMpdf}). 
The stability index of the attracting distribution is determined from the exponent of the tail behaviour of $p_s$  in (\ref{eqn:tail_behaviour})
\begin{equation}
\alpha^{*} = 2\nu + 1 = -2L/{E^{2}}.\label{eqn:alpha_conv}
\end{equation}
  For $0 < \nu < \frac{1}{2}$, $p_s$ has infinite variance and mean zero as discussed in Section \ref{subsec:properties_of_CAM}.
 The skewness parameter $\beta^*$ and scale parameter $\sigma^*$ are determined following \cite{Kuske2001}. First, $\beta^{*} = (h^{+} - h^{-})/({h^{+} + h^{-}})$ where $h^{\pm} = h(\pm 1) $ with $h(s)$ in (\ref{eqn:tail_behaviour}), which can be written simply
 as 
\begin{align}
\beta^{*} = \tanh\left( \frac{\pi g\nu}{b}\right). \label{eqn:beta_conv}
\end{align}
The scale parameter $\sigma^*$ is determined by comparing the arguments of the characteristic function of the normalized sum of $\{R_{j}\}_{j = 1}^{n}$ (\ref{eqn:CAM-to-stable}) in the limit $n \to \infty$ to those of an $\alpha$-stable random variable (\ref{eqn:stable_characteristic_function}) and is given by
\begin{align}
\sigma^{*} &= \left(\frac{(h^{+} + h^{-})\Gamma(1 - \alpha^{*})}{\alpha^{*}}\cos\left(\frac{\pi\alpha^{*}}{2}\right)\right)^{1/\alpha^{*}} \nonumber \\ &= \left(\frac{2\cosh(\pi g\nu /b )}{E^{(\alpha^{*} + 1)}\alpha^{*}\mathcal{N}}\Gamma(1 - \alpha^{*})\cos\left(\frac{\pi\alpha^{*}}{2}\right)\right)^{1/\alpha^{*}}. \label{eqn:sigma_conv}
\end{align}
 The analysis in \cite{Kuske2001} of the sum (\ref{eqn:CAM-to-stable}) gives an approximate expression for its PDF in the form of the attracting $\alpha$-stable PDF plus a correction of $O(n^{1 - 2/\alpha^*})$ for $\alpha^* \in (1, 2)$ and large $n$. Then $n^{1-2/\alpha^*}$ determines the rate of convergence to the $\alpha$-stable density for large $n$, with slower convergence for $\alpha^*$ closer to 2.
  If $\alpha^* = 2$ ($\nu = 1/2$), then the sum of $\{R_{j}\}_{j = 1}^{n}$ converges to a Gaussian distribution with a correction that is asymptotically $O(1/\log(n))$ for large $n$ \cite{Kuske2001}. We do not treat the case $\alpha^* = 2$ here, due to its logarithmic rate of convergence and the fact that the limiting distribution does not have power law tails.
  
The stochastic averaging approximation to which we now turn involves integrals of serially-dependent stochastic processes rather than sums of independent, identically-distributed random variables.  Our central ansatz is that Eqns. (\ref{eqn:alpha_conv})-(\ref{eqn:beta_conv}) determine the stability index and skewness parameter of such integrals, so that only the scale parameter needs to be computed.

\section{Approximation of a linear CAM process by an OULP}
\label{sec:reductionCAMtostable}

As $y_{t}$ is not an $\alpha$-stable process, we cannot directly apply the averaging results of \cite{Thompson2014} to obtain the weak approximation $X_{t}$ to $x_{t}$ in (\ref{eqn:slowSDE}).  However, we note that the solution $x_{t}$ of  (\ref{eqn:slowSDE}) for $t\in(0,T)$ involves the integral of $y_{t/\epsilon}$ over this interval. The connection between $\alpha$-stable random variables and 
sums of independent random variables with the stationary distribution of the CAM process $y_{t}$ described in Section \ref{sec:CAMandStable} suggests that for $\epsilon \ll 1$, the integral   $\int_{0}^{T}y_{s/\epsilon} ds$ should have a distribution that is close to $\alpha$-stable.  In this limit, we can determine the parameters of an OULP $z_{t}$ such that the distribution of $\int_{0}^{T}z_{s/\epsilon}\,ds$ approximates that of $\int_{0}^{T}y_{t/\epsilon}ds$, and then apply the stochastic averaging results from \cite{Thompson2014} to the two-timescale system $(x_{t},z_{t/\epsilon})$.  Our hypothesis is that the resulting stochastic process is a weak approximation to the slow variable $x_{t}$ in the original system $(x_{t},y_{t/\epsilon})$.

We cannot simply invoke the GCLT to make this claim however, since the integral of $y_{t}$ is not a sum of independent, identically-distributed random variables.  Demonstrating this result involves arguing that the integral $\int_{0}^{T}y_{s/\epsilon}ds$ can be decomposed into the sum of a large number of effectively independent and identically-distributed random variables.   As well, previous stochastic averaging results for fast-slow systems such as \cite{Mona11} and \cite{Thompson2014} indicate that details of the serial dependence of the fast process appear directly in the stochastic averaging approximation.  Specifically, 
when an $\alpha$-stable process with characteristic decay time scale $\tau$ drives the fast process, then the $\alpha$-stable process in the the stochastic averaging approximation depends 
%on the characteristic decay time scale from the serial dependence of the fast process. 
on $\tau$.  This dependence appears in the scale parameter used
in the approximation, but not in its stability index and skewness parameter.
 Thus, we expect that any $\alpha$-stable stochastic forcing replacing the fast CAM noise process would depend on the memory of $y_{t}$ through the scale parameter.
 
 With these facts in mind, we proceed to investigate the hypothesis that the system (\ref{eqn:slowSDE}), (\ref{eqn:CAMsde}) can be reduced by relating the fast driving linear CAM process $y_{t/\epsilon}$ to  an OULP, $z_{t/\epsilon}$ with approximately equally-distributed integrals over the interval $(0,T)$.  Then the stochastic averaging result from \cite{Thompson2014} is used to approximate the weak properties of the fast fluctuations in the equation for  $x_{t}$ when there is a separation of the  slow and fast time scales for $\epsilon \ll 1$. 
%Since it is the properties of the integral of $y_{t/\epsilon}$ that appear in the equation for $x_{R}$, we find conditions under which it can be approximated by the integral of an OULP  $z_{t/\epsilon}$. 
 In order to identify the appropriate parameter values in the SDE for
$z_{t/\epsilon}$, we review the statistical properties of $\int_{0}^{t}z_{s/\epsilon}\,ds$ as derived in \cite{Thompson2014}.  These are complemented
by the conditions under which the integral  $\int_{0}^{t}y_{s/\epsilon}\,ds$, can be weakly approximated by an $\alpha$-stable random variable.

%We conclude this section by using the (N+) stochastic averaging approximation (\ref{eqn:N+approx_general}) to complete our stochastic averaging approximation for $x_{t}$.

\subsection{The distribution for the integral of an OULP}
\label{subsec:OULP_approx_for_y}

We derive the characteristic function for the integral $v_{t} = \int_{0}^{t}\,z_{s/\epsilon}\,ds$
 ($dv_{t} = z_{t/\epsilon}\,dt$) 
 where $z(t)$ is the OULP in (\ref{eqn:zt_dynamics}).
 The $(v_{t},z_{t/\epsilon})$ system has a joint time-dependent PDF $P(v,z,t)$  satisfying the FKE,
\begin{eqnarray}
\frac{\dd P}{\dd t} &=& \left(\frac{\theta}{\epsilon}\frac{\dd}{\dd z} -\frac{\dd}{\dd v}\right)(zP) + \frac{\sigma_z^{\alpha}}{\epsilon}\mathcal{D}_{z}^{(\alpha,\beta)}P, \label{eqn:FKE_vz}\\
\mathcal{D}_{z}^{(\alpha,\beta)} &=& -\frac{\beta}{\cos(\pi\alpha/2)} \frac{\partial^{\alpha}}{\partial y^{\alpha}} +(1-\beta)\frac{\partial^{\alpha}}{\partial y^{\alpha}}.
\end{eqnarray}
where $P(v,z,0) = \delta(v)\delta(z - z_{0}) ,\,t > 0$. The operator $\mathcal{D}_{z}^{(\alpha,\beta)}$ is the fractional differentiation operator \cite{Chaves1998}.
% that is defined in terms of a Fourier transform with respect to $z$, corresponding to $\alpha$-stable noise driving with stability index $\alpha$ and skewness parameter $\beta$ \cite{Chaves1998}. 
As shown in Appendix \ref{calc_details}, we obtain the joint characteristic function $\psi_{v,z}(m,k,t) $ for $(v_t, z_{t/\epsilon})$ from 
 (\ref{eqn:FT_FKE_vz}), the Fourier transform of (\ref{eqn:FKE_vz}).  We find that $\psi_{v,z}(m,k,t) = \psi_{z}(k)\psi_{v}(m,t)$  for $t=O(1)$ and $\epsilon \ll 1$, where
 \begin{eqnarray}
\psi_{z}(k) & =& \exp\left(\frac{\sigma_z^{\alpha}}{\alpha\theta}|k|^{\alpha}\Xi(k;\alpha,\beta)\right) + O(\epsilon) \label{eqn:char_func_z}\\
\psi_{v}(m,t) &=& \exp\left[i\frac{\epsilon m z_{0}}{\theta}(1 - e^{-\theta t/\epsilon})  - \frac{\epsilon^{\alpha - 1}\sigma_z^{\alpha}t}{\theta^{\alpha}}|m|^{\alpha}\Xi(m;\alpha,\beta)\right] +O(\epsilon) \label{eqn:char_func_v}
\end{eqnarray}
for $\Xi(m;\alpha,\beta)$ as given in (\ref{eqn:Xi}). The expression for $\psi_z$ in (\ref{eqn:char_func_z}) is the characteristic function corresponding to the stationary density of the OULP $z_t$. We note that the ratio $\sigma_z/ \theta^{1/\alpha}$ appears in $\psi_z(k)$, analogous to the expression for the standard deviation in a stationary OU process.  For $z_{0}=0$,
% assuming it is drawn from the stationary distribution of $z_{t}$ with zero mean,
 the form of $\psi_v$ is the same as the characteristic function of an $\alpha$-stable random variable with mean $0$ and the scale parameter  $\epsilon^{1-1/\alpha}\sigma_z\theta^{-1}t^{1/\alpha}$. As convergence in characteristic function implies convergence in distribution, the asymptotic distribution for $v_{t}$ is $\alpha$-stable, 
\begin{equation}
v_{t} = \int_{0}^{t}z_{s/\epsilon}\,ds \quad \underset{D}{\to} \quad  \mathcal{S}_{\alpha}\left(\beta,\epsilon^{\gamma}t^{1/\alpha}\frac{\sigma_z}{\theta}\right) \label{eqn:integral_of_z},
\end{equation}
where $\gamma = 1 - 1/\alpha$.

\subsection{The distribution of the integral of CAM noise}
\label{subsec:CAM2stablenoise}
\begin{figure*}
\includegraphics[width=0.49\textwidth]{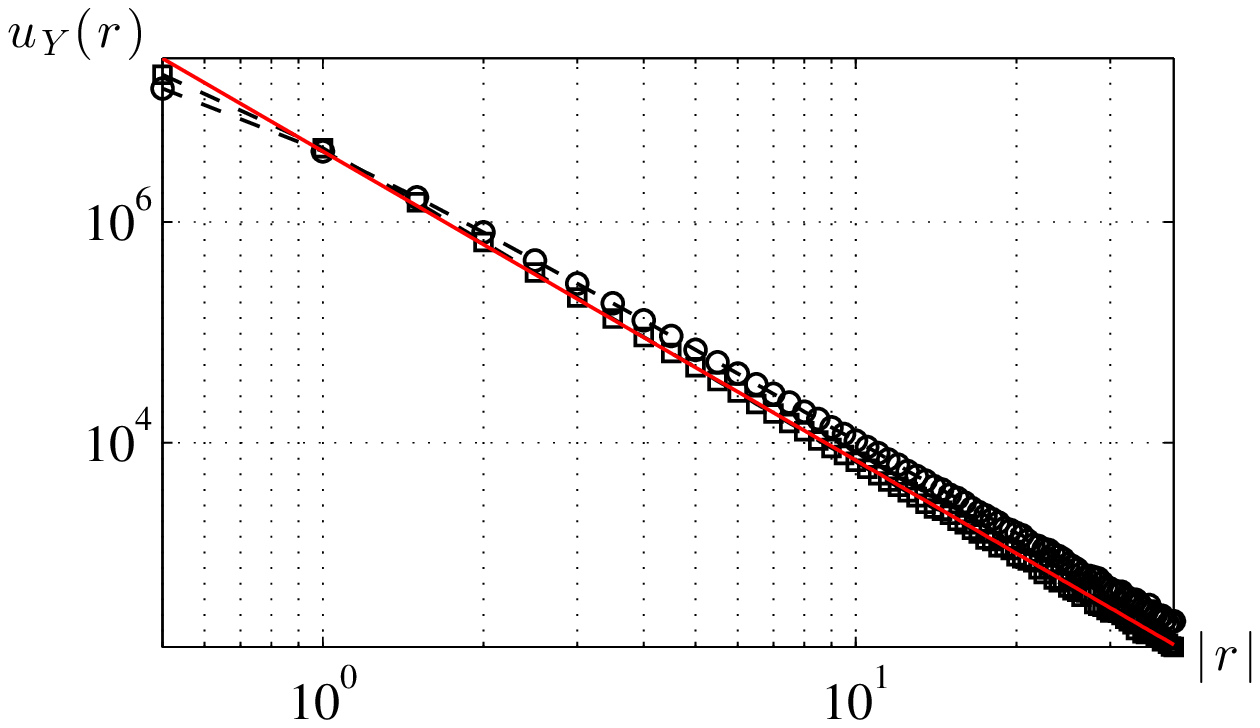}
\includegraphics[width=0.49\textwidth]{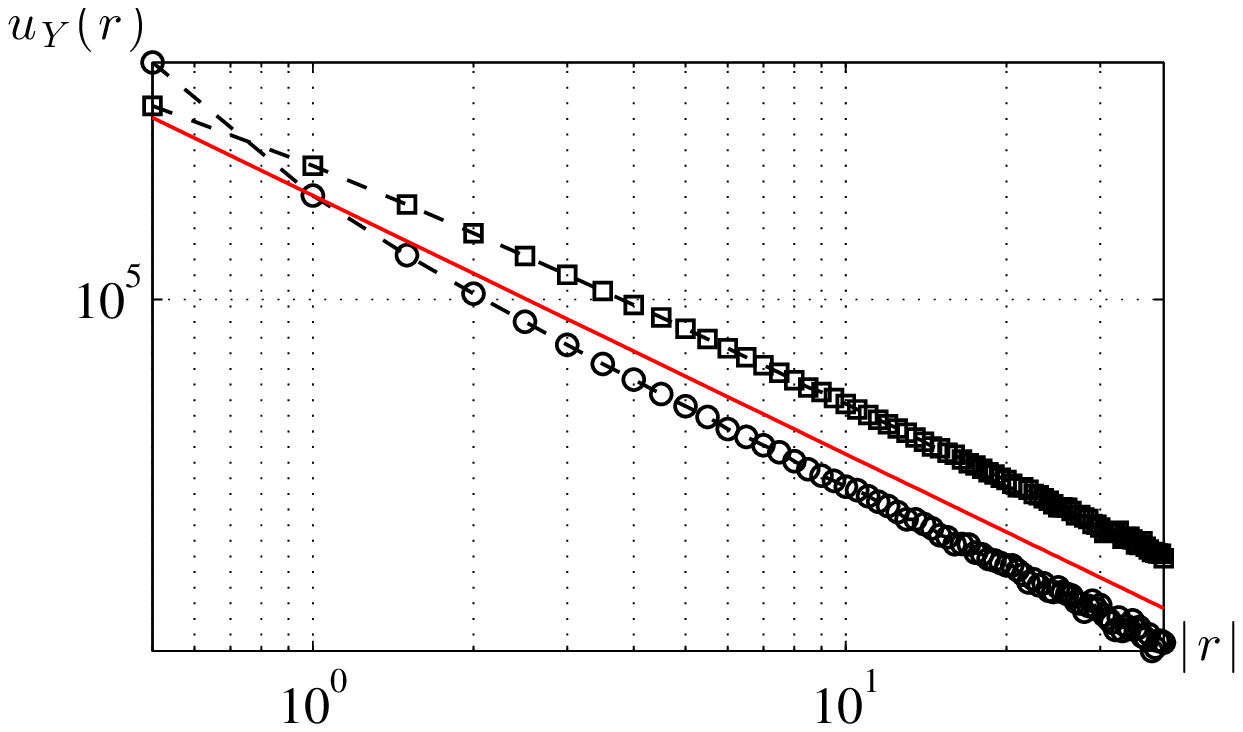}\\
\includegraphics[width=0.49\textwidth]{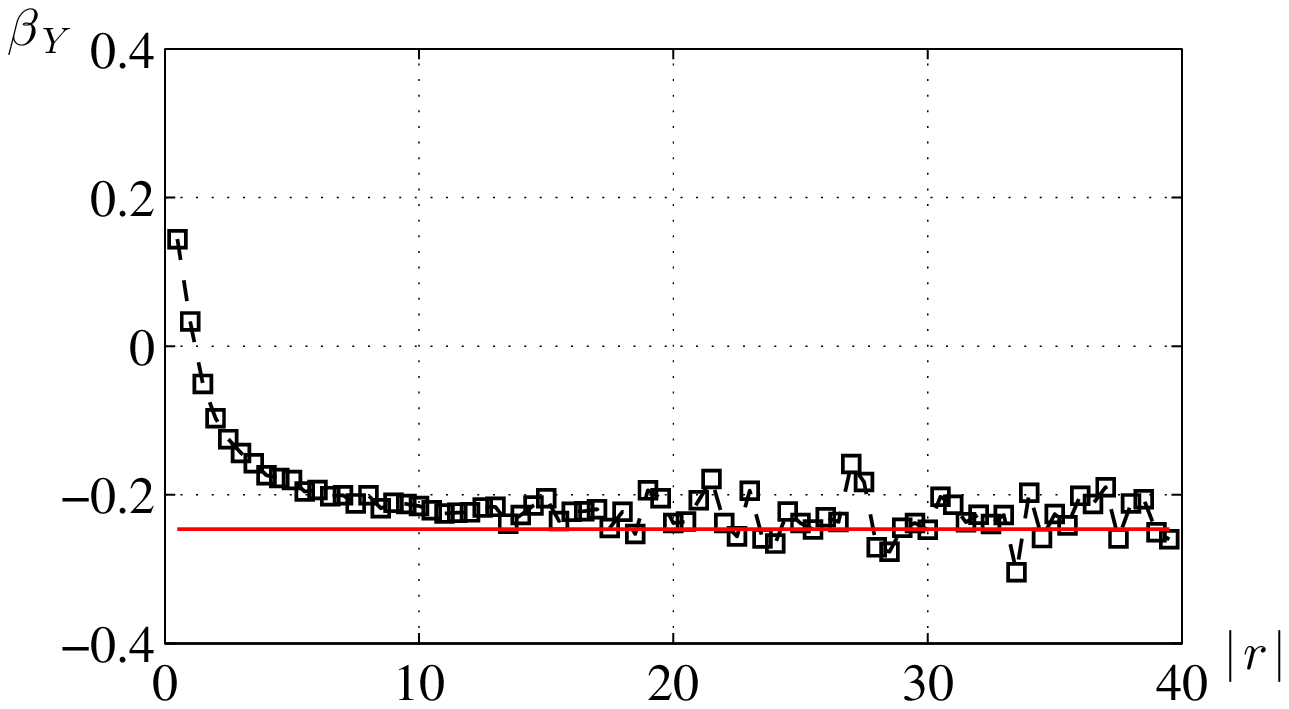}
\includegraphics[width=0.49\textwidth]{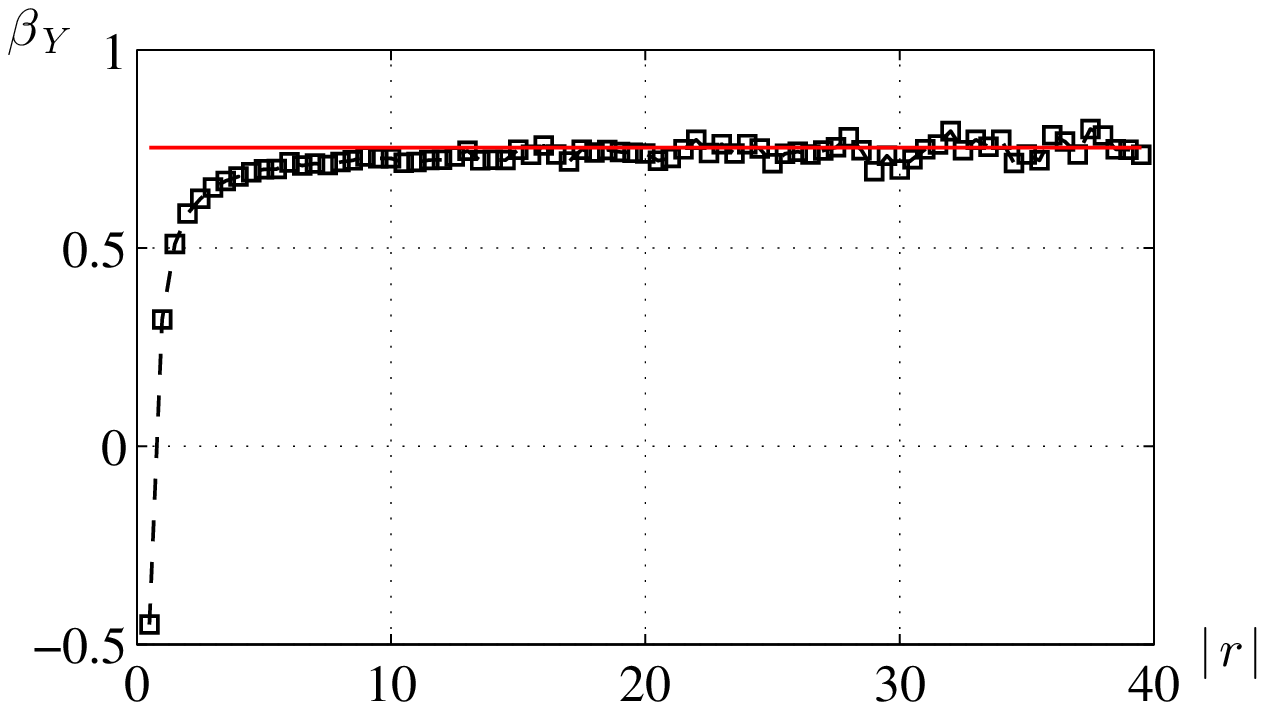}
\caption{Top: Histogram-based comparison of tail behaviour of $u_Y(r)$ (squares for $r$ positive, circles for $r$ negative) to the plot of $cr^{-\alpha^*-1}$, $c$ constant (red line). Bottom: The corresponding plots of $(u_Y(r) - u_Y(-r))/(u_Y(r) + u_Y(-r)) \sim \beta^*$ from (\ref{eqn:qval}) for $|r|$ large. The red line shows the value of $\beta^*$ corresponding to the following specific parameters --  Left: $\alpha^* = 1.5, g = 0.25, b = 0.2, L = -1$; Right: $\alpha = 1.8, g = -0.1, b = 0.5, L = -1$. Estimates are based on  $10^8$ realizations of $Y_j$ (\ref{eqn:integral_yseps}), with $\Delta = \epsilon = 1$.}\label{fig:beta}
\end{figure*}

We begin by considering the integral $\int_{0}^{t} y_{s/\epsilon}\,ds$ as a sum of $N_{Y}$ partitions of length $\Delta$, denoted by $\{Y_{j}\}_{j = 1}^{N_{Y}}$,
\begin{eqnarray}
&\int_{0}^{T}&y_{s/\epsilon}\,ds = \sum_{j = 1}^{N_{Y}}Y_{j}, \label{eqn:integral_yseps} \\ &Y_{j}& = \int_{(j - 1)\Delta}^{j\Delta}y_{s/\epsilon}\,ds = \epsilon\left(\int_{(j-1)\Delta/\epsilon}^{j\Delta/\epsilon}y_{\hat{s}}\,d\hat{s}\right), \label{eqn:integral_yseps}
\end{eqnarray}
where $T = N_{Y}\Delta$. 
We make the following ansatz for the tail behaviour of $Y_j$: {\sl for any fixed non-zero value $\Delta$ and sufficiently small $\epsilon$,
$Y_j$ has a PDF $u_{Y}$, where
\begin{equation}
u_{Y}(r) \sim \begin{cases} q^{-}|r|^{-(\alpha^{*} + 1)} &{\mbox{for $r < -a$ }} \\ q^{+}r^{-(\alpha^{*} + 1)} &\mbox{for $r > a$} \end{cases}, \label{eqn:dist_RR}
\end{equation}
for some $ a > 0$ and $q^{+},q^{-}>0$.}. The ansatz for $u_Y$ follows from noting that the  integral (\ref{eqn:integral_yseps}) for $Y_j$ can be approximated by a sum of random variables 
$y_m$,  for $Mh=\Delta$ and sufficiently small $h$. In Appendix \ref{app:uY} we use a Euler-Maruyama approximation of $y_m$ in (\ref{eqn:CAMsde}),  with $dW_{i,t}$ approximated by sequences of independent random variables $\xi_{i,m}\sim N(0,\sqrt{h})$  for $i=1,2$. We express $u_Y(Y_j)$ in terms of $p_s(y)$ via a linear relationship between $Y_j$ and $y_1$ (\ref{Yjtoy1}).
  This linear transformation leads not only to the exponent $\alpha^*$ in the ansatz (\ref{eqn:dist_RR}), but also to  an approximation of 
a parameter $\beta_Y$ by the skewness parameter $\beta^*$ in (\ref{eqn:tail_behaviour}). Specifically, if $q^{\pm}$ are constant multiples of $h^{\pm}$ in (\ref{eqn:tail_behaviour}), then
\begin{eqnarray}
\beta_Y \equiv \frac{q^{+} - q^{-}}{q^{+} + q^{-}} \approx \frac{h^{+} - h^{-}}{h^{+} + h^{-}} = \beta^{*} \, .\label{eqn:qval}
\end{eqnarray}
Figure \ref{fig:beta} illustrates the ansatz (\ref{eqn:dist_RR}) and the skewness parameter $\beta_Y \approx \beta^*$ for representative examples of $Y_j$. The results show that indeed the values of $Y_{j}$ have the tail behaviour assumed by the ansatz that we assume. The rate of decay of both the positive and negative tails is consistent with a power law with exponent $-1-\alpha^{*}$ and the relative weight of the tails is as predicted for $|r|$ sufficiently large. Thus, we have empirical evidence for the validity of our ansatz, ((\ref{eqn:dist_RR}) and (\ref{eqn:qval})).   Note that while we can predict the value of the ratio $(q^{+}-q^{-})/(q^{+}+q^{-})$ from the value of $\beta^{*}$, we do not have expressions for the values of the numerator or denominator separately.

The GCLT implies that a sum of independent variables $Y_{j}$ distributed according to (\ref{eqn:dist_RR}) converges in distribution to an $\alpha$-stable random variable with stability index $\alpha^{*}$ and skewness parameter $\beta_{Y}\approx \beta^*$ as $N_Y \to \infty$.
Thus, we expect to be able to approximate the distribution of $\int_{0}^{T}y_{s/\epsilon}\,ds$ with that of an
$\alpha$-stable random variable with an error that decreases with $N_Y$, if the following two conditions are satisfied.
\begin{itemize}
\item{ Condition A) { \sl The value of $\Delta$ must be large enough such that $Y_{j}$ and $Y_{j+1}$ are (to a good approximation) independent for all $j \in \{1,2,\dots, N_{Y}-1\}$.}}
\item{Condition B){ \sl The number of partitions $N_{Y} = T/\Delta$ must be large enough such that the distribution of the normalized sum of $\left\{Y_{j}\right\}_{j = 1}^{N_{Y}}$ is well-approximated by that of an $\alpha$-stable random variable.}}
\end{itemize}
If Conditions A) and B) are satisfied for $T = N_{Y}\Delta$, then $N_Y^{-1/\alpha^*}\sum_{j = 1}^{N_{Y}}Y_{j}$ can be approximated by a $\alpha$-stable random variable  (\ref{eqn:char_funcS}).
Furthermore, Conditions A) and B) taken by themselves appear to suggest that $T$ should be large; however, when considered within the larger context of approximating $x$ in (\ref{eqn:slowSDE}) on the slow time scale, we see that $T$ must be of the same order of the characteristic time scale of $x$.  In fact, these conditions can be satisfied by taking $\epsilon$ to be sufficiently small. The asymptotic approximations that use Conditions  A) and B) are discussed below.

\begin{enumerate}[A)]
\item {\sl The value of $\Delta$ must be large enough such that $Y_{j}$ and $Y_{j+1}$ are (effectively) independent for all $j \in \{1,2,\dots, N_{Y}-1\}$.}

For the GCLT to apply to the sum $\sum_jY_{j}$,  the individual terms $Y_j$ must be independent.  We intuitively expect $Y_j$ and $Y_{j+1}$ to be asymptotically independent (and thus also $Y_j$ and $Y_{j+m}$ for $m>1$) provided that $\Delta/\epsilon \gg \tau_{y}$, for $\tau_{y}$   a characteristic memory time scale of $y_{t}$ and $\Delta/\epsilon$ in the definition of $Y_j$ (\ref{eqn:integral_yseps}). For example, $\tau_y$ could be the $e$-folding time of the ACD function: $ACD_{y}(\tau_{y}) = e^{-1}ACD_{y}(0)$.  For $\operatorname{ACD}_y(\tau)$ decreasing with $\tau$ as shown in Figure \ref{fig:sample_acds}, we see that $\tau_y= O(1)$. 

 Two random variables with infinite variance $Y_k$ and $Y_l$ are independent if both the codifference ($CD$) and the cosum ($CS$)  vanish \cite{Gross1994, Wylomanska2015}, where
\begin{eqnarray}
 CD(Y_k,Y_l) &=& \log(\expected{\exp(iY_k-iY_l)}) - \log(\expected{\exp(iY_k)}) - \log(\expected{\exp(-iY_l)}), \nonumber\\
 CS(Y_k,Y_l) &=& CD(Y_k,-Y_l)\, .
\end{eqnarray}
 Hence, we consider the asymptotic behaviour of $CD(Y_j,Y_{j+1})$ and $CS(Y_j,Y_{j+1})$, which decrease with increasing $\Delta/\epsilon$ as illustrated in Figure \ref{fig:codiff} for example CAM processes. It should be noted that the choice to normalize the values of $Y_{j}$ in the figure is done for stylistic purposes and does not affect the relative values of the cosum and codifference functions for different values of $\Delta/\epsilon$. We expect that the codifference for an $\alpha$-stable random variable with itself is equal to 2 when the scale parameter is normalized to 1 \cite{Wylomanska2015}, as the results in Figure \ref{fig:codiff}. For $CD(Y_j,Y_{j+1})$ and $CS(Y_j,Y_{j+1})$ approaching zero, we treat $Y_{k}$ for $k>j+1$ as independent $Y_j$, treating all $Y_j$ effectively independent. It is clear that for $|L| = O(1)$, we can expect effective independence of $Y_{j}$ and $Y_{j+1}$ provided that $\Delta/\epsilon\gtrsim 5$.

\begin{figure}
\includegraphics[width=0.49\textwidth]{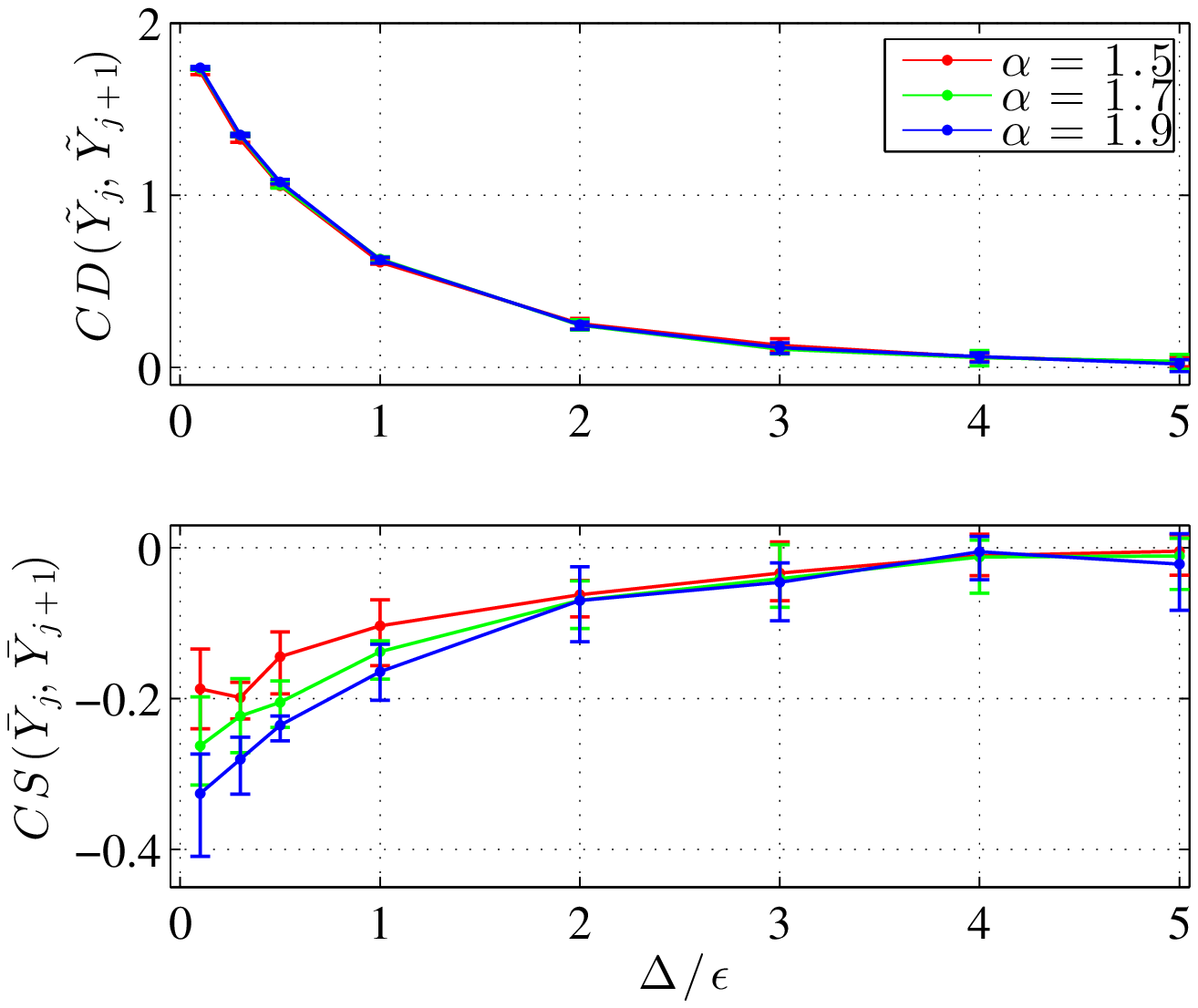}
\includegraphics[width=0.49\textwidth]{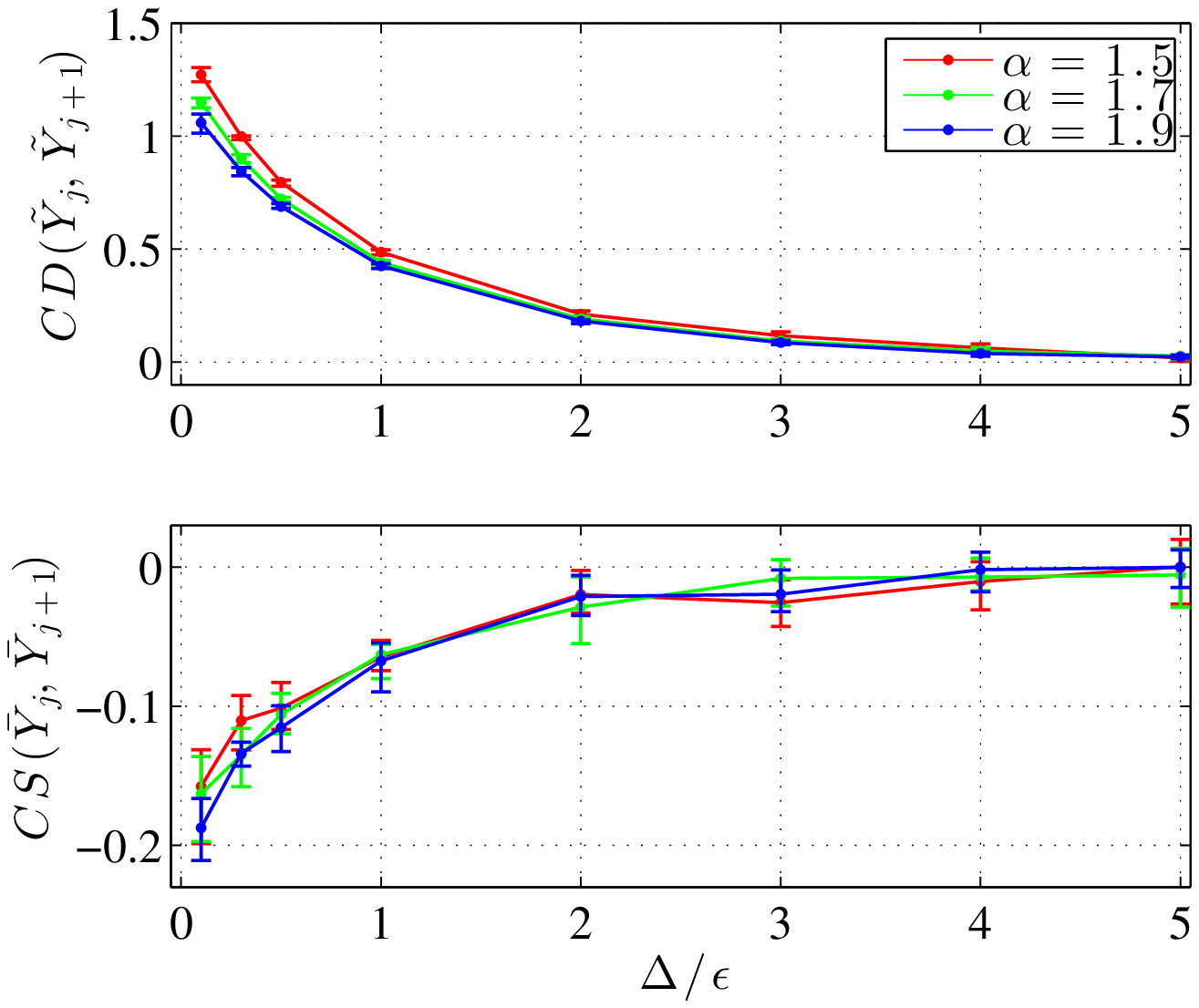}
\caption{Plots of the codifference $CD$ and cosum $CS$ of ${Y}_{j}/s$ and $Y_{j+1}/s$ as a function of the length of time integral $\Delta/\epsilon$ as given in (\ref{eqn:integral_yseps}) where $s$ is the numerically estimated scale parameter of $Y_{j}$. For both plots, $L = -1$, and $\alpha$ is as indicated in the legend. Twenty series of $10^{4}$ pairs of $(Y_{j},Y_{j+1})$ are sampled for each value of $\alpha$ considered in each figure. The ends of the error bars denote the 25th and 75th percentiles. {Left:} $g = 0.1,\,b = 1$.  {Right:} $g = 0.4,\,b = 0.6$.}
\label{fig:codiff}
\end{figure}

\item {\sl  The number of partitions $N_{Y} = T/\Delta$ must be large enough such that the distribution of the sum of $\left\{Y_{j}\right\}_{j = 1}^{N_{Y}}$ is well-approximated by that of an $\alpha$-stable random variable.}
To characterize the error in relating  $\sum _{j = 1}^{N_Y}Y_{j}$ to an $\alpha$-stable random variable, we review the results from \cite{Kuske2001} that gives the distribution of a sum $S_{N_Y}$ of heavy-tailed random variables $Y_j$, 
\begin{eqnarray}
S_{N_{Y}} = \left(N_{Y}^{-1/\alpha^{*}}\sum_{j = 1}^{N_{Y}}Y_{j}\right) \, ,\label{SNdef}
\end{eqnarray}
via its characteristic function $\psi_S$. The rate of convergence to this distribution depends on its tail behaviour. For $\Delta/\epsilon$ large enough for effective independence of subsequent values of $Y_{j}$, we can write the PDF $p_S(r)$ of $S_{N_{Y}}$  as
\begin{equation}
p_S(r) = \frac{1}{2\pi}\int_{\mathbb{R}}\exp\left(-ikr\right)\psi_{S}\left(k\right)\,dk, \quad \mbox{where}\quad \psi_{S}(k) \sim \psi_{Y}^{N_Y}\left(\frac{k}{N_{Y}^{1/\alpha^{*}}}\right). \label{eqn:Yi_sum_dist_fourier}
\end{equation}
and  $\psi_{Y}$ is the characteristic function of $Y_{j}$.  Then for large $N_Y$,  using the result of \cite{Kuske2001} for a sum of independent, identically distributed heavy-tailed random variables, the characteristic function $\psi_S(k)$ can be approximated from an expansion of $\psi_Y$ for small argument.
Using the expansion (\ref{eqn:apprx_charfuncY}) for $\psi_Y$ as derived in Appendix  \ref{subsec:asymp_charfunc} yields the approximation
\begin{equation}
\psi_{S}(k) = \exp\left[ -\left(\frac{(q^{+} + q^{-})\Gamma(1 - \alpha^{*})}{\alpha^{*}}\cos\left(\frac{\pi\alpha^{*}}{2}\right)\right)|k|^{\alpha^{*}}\Xi\left(k;\alpha^{*},\beta^{*}\right) - \frac{Q}{N_{Y}^{2/\alpha^{*} - 1}}k^{2} + O\left(\frac{k^{3}}{N_{Y}^{3/\alpha^{*} - 1}}\right)\right] \label{psi_Sexp}
\end{equation}
for $Q$ a constant. As $N_Y$ increases, (\ref{psi_Sexp}) approaches the form of a characteristic function for an $\alpha$-stable distribution as given
in (\ref{eqn:stable_characteristic_function})
\begin{widetext}
\begin{equation}
\psi_{S}(k) \to \exp\left[ -\left(\frac{(q^{+} + q^{-})\Gamma(1 - \alpha^{*})}{\alpha^{*}}\cos\left(\frac{\pi\alpha^{*}}{2}\right)\right)|k|^{\alpha^{*}}\Xi\left(k;\alpha^{*},\beta^{*}\right)\right] \quad \mbox{as $N_{Y} \to \infty$.} \label{eqn:char_funcS}
\end{equation}
\end{widetext}
Then  the distribution of $S_{N_{Y}}$ converges to an $\alpha$-stable distribution, implied by the pointwise convergence of their characteristic functions (as per L\'{e}vy's continuity theorem \cite{Feller1966II}).  The error term in (\ref{psi_Sexp}) has coefficient $N_{Y}^{1-2/\alpha^*}$, characterizing the rate of convergence of $S_{N_{Y}}$ to an $\alpha$-stable distribution, as noted in Section \ref{sec:CAMandStable}.   This error increases substantially for $\alpha^* <2$ and $(2-\alpha^*) \ll 1$. 
\end{enumerate}
Putting these results together, it follows from (\ref{eqn:char_funcS}) that we can specify the parameters for the approximate distribution of $S_{N_Y}$ as
\begin{eqnarray}
&S_{N_{Y}}&\,\underset{D}{\to}\,\mathcal{S}_{\alpha^{*}}(\beta^{*},\sigma_Y), \label{sigma_Ydef} \\ 
&{\sigma_Y}& = \left(\frac{(q^{+} + q^{-})\Gamma(1 - \alpha^{*})}{\alpha^{*}}\cos\left(\frac{\pi\alpha^{*}}{2}\right)\right)^{1/\alpha^{*}}, \label{eqnarray}
\end{eqnarray}
and $T = N_{Y}\Delta$, and $N_{Y}, \Delta$ are sufficiently large for Conditions A) and B) to hold.   Note that because we do not know the value of $q^{+}+q^{-}$, Eqn. (\ref{sigma_Ydef}) does not allow us to compute the value of $\sigma_{Y}$.

\begin{figure*}[t]
\includegraphics[width=0.98\textwidth]{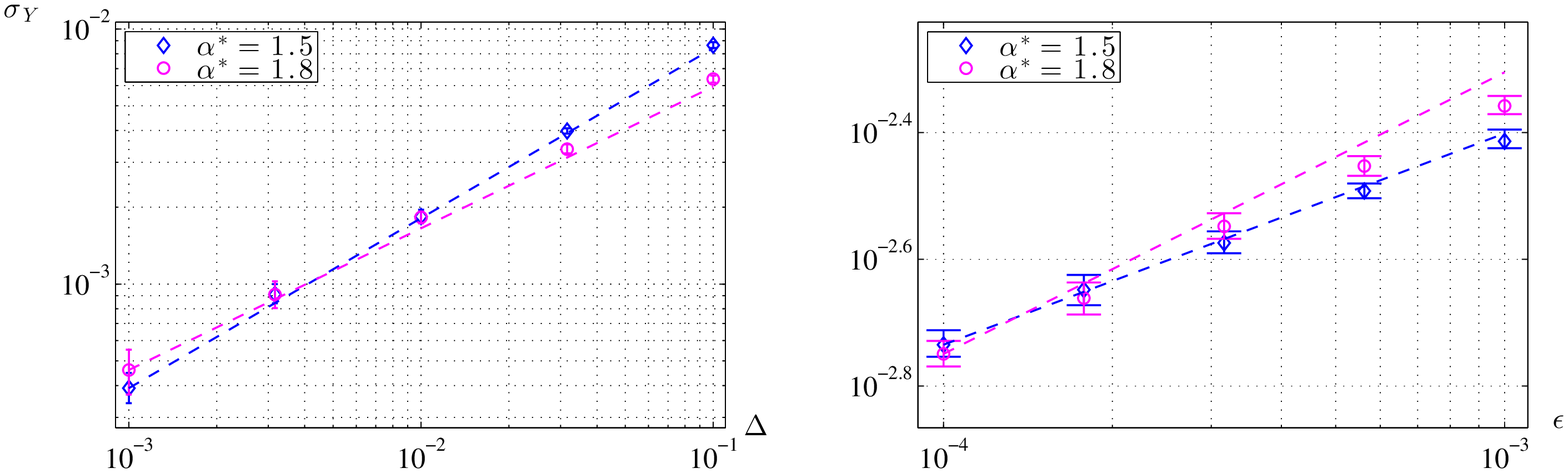}
\caption{
Numerically estimated values of the scale parameter $\sigma_Y$ (open circles), based on realizations of $S_{j} = \int_{0}^{T}y_{s/\epsilon}\,ds$, as functions of $\epsilon$ and $\Delta$ for different values of $\alpha^{*}$. The fixed parameters are $(L,g,b) = (-1,0.1,0.5)$ and $N_{Y} = 100$. {Left}: Dashed lines are proportional to $\Delta^{1/\alpha^{*}}$ and $\epsilon = 10^{-4}$. {Right}: dashed lines are proportional to $\epsilon^{\gamma^{*}}$ and $T = 1$. The error bars indicate the upper and lower quartiles of the estimates of $\sigma_{Y}$. Each bar chart is based on 30 estimates of $\sigma_{Y}$, each of which is generated from $N_{S} = 1000$ simulated instances of $S_{j}$.}
\label{fig:sigma_dependence}
\end{figure*}

Now, we demonstrate that $\sigma_Y$ is proportional to $\epsilon^{\gamma^{*}}\Delta^{1/\alpha^{*}}$ for 
$\gamma^* = 1-1/\alpha^*$, which allows us to determine the dependence on $\epsilon$ and $T$ of the scale parameter of $\int_{0}^{T}y_{s/\epsilon}\,ds$.
The integral (\ref{eqn:integral_yseps}) can be expressed by multiplying the sum $S_{N_{Y}}$ by $N_{Y}^{1/\alpha^{*}}$ which gives us an $\alpha$-stable distribution for $\int_{0}^{T}y_{s/\epsilon}\,ds$:
\begin{equation}
 \int_{0}^{T}y_{s/\epsilon}\,ds \simeq N_{Y}^{1/\alpha^{*}}S_{N_{Y}}\,\underset{D}{\to}\,\mathcal{S}_{\alpha^{*}}(\beta^{*},N_{Y}^{1/\alpha^{*}}\sigma_Y) = \mathcal{S}_{\alpha^{*}}\left(\beta^{*},\left(\frac{T}{\Delta}\right)^{1/\alpha^*}\sigma_Y\right) \label{eqn:int_yeps_intermed1}
\end{equation}
It follows that
\begin{equation}
\int_{0}^{T}y_{\hat{s}}\,d\hat{s}\,\underset{D}{\to}\,\mathcal{S}_{\alpha^{*}}\left(\beta^{*},\epsilon^{-\gamma^*} \left(\frac{T}{\Delta}\right)^{1/\alpha^*}{\sigma_Y}\right),\label{eqn:int_y_no_eps_dep}
\end{equation}
for sufficiently large $T$. The integral $\int_{0}^{T}y_{\hat{s}}\,d\hat{s}$ has no dependence on $\epsilon$ or $\Delta$, so that
\begin{equation}
\sigma_Y \propto\epsilon^{\gamma^{*}}\Delta^{1/\alpha^*}. \label{sigma_Yprop}
 \end{equation}
This dependence of ${\sigma}_{Y}$ on $\epsilon$ and $\Delta$ is illustrated in Figure \ref{fig:sigma_dependence}.
Using this result, it is  useful to rewrite (\ref{eqn:int_yeps_intermed1}) as
\begin{equation}
\int_{0}^{T}y_{s/\epsilon}\,ds \,\underset{D}{\to}\, \mathcal{S}_{\alpha^{*}}\left(\beta^*,\epsilon^{\gamma^{*}}{T}^{1/\alpha^*}\Sigma\right), \quad \Sigma =\frac{{\sigma_Y}}{\epsilon^{\gamma^{*}}\Delta^{1/\alpha^{*}}}. \label{eqn:int_y_to_stable}
\end{equation}
Note that by (\ref{sigma_Yprop}),  $\Sigma$ is independent of $\epsilon$ and $\Delta$ and depends only on the parameters of the CAM noise in (\ref{eqn:CAMsde}). However,  (\ref{sigma_Ydef}) and (\ref{sigma_Yprop}) do not provide the value of $\Sigma$, since the coefficients $q^{\pm}$ are not specified.
We numerically estimate $\Sigma$ via the empirical characteristic function obtained from simulations of $Y_j$, as described in Appendix \ref{sec:estim_sigma}. 

Finally, we obtain the approximation result that on a sufficiently long time  $T$, 
\begin{equation}
\int_{0}^{T}y_{s/\epsilon}\,ds \simeq \,\int_{0}^{T}z_{s/\epsilon}\,ds \quad \underset{D}{\to} \quad \mathcal{S}_{\alpha^{*}}\left(\beta^{*},\epsilon^{\gamma^{*}} T^{1/\alpha^{*}}\Sigma\right) \, ,\label{etaz_equiv}
\end{equation}
if the
parameters  in (\ref{eqn:zt_dynamics}) are taken to be
\begin{equation}
dL_{t}^{(\alpha,\beta)}=dL_{t}^{(\alpha^*,\beta^*)}, \qquad \sigma_z = \Sigma\theta \, .\label{eqn:theta_true}
\end{equation}

Here $\alpha^{*}$ and $\beta^{*}$ are given by Eqns (\ref{eqn:alpha_conv})-(\ref{eqn:beta_conv}) and $\Sigma$ is estimated empirically as described in Appendix \ref{sec:estim_sigma}.   Note that for the averaging approximation the choice of $\theta$ is arbitrary, as only $\Sigma$ enters the expression for the scale parameter of the integral of $y_{t/\epsilon}$.  A reasonable choice is to use a characteristic inverse timescale of Eqn. (\ref{eqn:CAMsde}) for $y_{t}$ and set  $\theta = L + E^{2}/2$.

\section{The stochastic averaging approximation}
\label{sec:stochastic_averaging_approx_final}

We now use the results from %Section \ref{subsec:OULP_approx_for_y}   
the previous section to find the SDE for a slow process $X(t)$ that weakly approximates $x(t)$.  To accomplish this
we  find a weak approximation for the fast fluctuations that involve $y_{t/\epsilon}$ in the equation for $x(t)$.  Following \cite{Thompson2014}, 
 we make a change of variable $\eta_{t} = \mathcal{U}(x_{t})$ for $\mathcal{U}'(x) = \frac{1}{f_{2}(x)}$, so that $y_{t/\epsilon}$ enters as an additive term in the equation for $\eta(t)$:
\begin{eqnarray}
d\eta_{t} &= \mathcal{U}'(x_{t})\,dx_{t} = \tilde{f}(\eta_{t})\,dt + \epsilon^{-\rho}y_{t/\epsilon}\,dt, \\ 
&\mbox{where} \quad \tilde{f}(\eta) = \frac{f_{1}(\mathcal{U}^{-1}(\eta))}{f_{2}(\mathcal{U}^{-1}(\eta))}.
\label{eqn:dynamics_eta}
\end{eqnarray}
Notice that the transformation $\mathcal{U}$ is invertible, since $f_{2}(x) \ne 0$ for any $x$ in the domain. We write the dynamics (\ref{eqn:dynamics_eta}) in integral form,
\begin{equation}
\eta_{t} = \eta_{0} + \int_{0}^{t}\tilde{f}(\eta_{s})\,ds + \epsilon^{-\rho}\int_{0}^{t}y_{s/\epsilon}\,ds,  \quad t > 0, \label{eqn:slowSDE_integral}
\end{equation}
and consider the integral of $y_{t/\epsilon}$ over a time interval of length $T$, $\int_{0}^{T}y_{s/\epsilon}\,ds$.

Then we approximate $\int_{0}^{T}y_{s/\epsilon}\,ds$ in the SDE for $\eta(t)$ (\ref{eqn:slowSDE_integral})  using $\int_{0}^{T}z_{s/\epsilon}\,ds$  where $z_{t}$ satisfies (\ref{eqn:zt_dynamics}) with $\theta$ and $\sigma_z$ in (\ref{eqn:theta_true}),
\begin{equation}
d\eta_{t} \approx \left(\tilde{f}(\eta_{t}) +  \epsilon^{-\rho} z_{t/\epsilon}\right)\,dt.\label{eqn:dynamics_eta_with_z}
\end{equation}
While $T$ must be sufficiently large for (\ref{etaz_equiv}) to hold,  we also require that the approximation is appropriate on the characteristic timescale $\tau_x$ of $x_{t}$, consistent with the overarching goal  to find an approximate equation for the dynamics on the slow variable $x_t$.  
%In Section \ref{discuss_num} we discuss how the requirement $T\approx \tau_x$ appears in the numerical results.

The system (\ref{eqn:dynamics_eta_with_z}) and (\ref{eqn:zt_dynamics}) is a slow-fast system for $\eta_{t}$ and $z_t$ for which we seek a stochastic averaging approximation, as an intermediate step to weakly approximating $x_t$.
For completeness, we give the key components of the (N+) approximation based on \cite{Thompson2014} for the general slow-fast system,
\begin{eqnarray}
du_{t}&=& f(u_{t},\hat{z}_{t})\,dt = \left(f_{1}(u_{t}) + \epsilon^{-\rho}f_{2}(u_{t})\hat{z}_{t}\right)\,dt, \label{eqn:u_dynamics_stable_averaging}\\
d\hat{z}_{t}&=& \epsilon^{-1/\alpha}g_{1}(u_{t})\,dt - \epsilon^{-1}g_{2}(u_{t})\hat{z}_{t}\,dt + \epsilon^{-1/\alpha}\sigma\,dL_{t}^{(\alpha,\beta)}
\label{eqn:hatz}
\end{eqnarray}
where $ t \ge 0, 0 < \epsilon \ll 1$, $\rho$ is a constant, $f_{2}({u}) \ne 0$ and $g_{2}({u}) > 0$ for any ${u}$ in the domain of $u_{t}$.  Note that  (\ref{eqn:zt_dynamics})
is a special case of (\ref{eqn:hatz}).
%The dynamics of $x_{t}$ are expressed in terms of a slow mean drift that is independent of $z_{t}$ and a fast perturbation, $v_{t} = v_{0} + \int_{0}^{t}f(x_{s},z_{s}) - \expectedcond{f(x_{s},z_{s})}{z_{s}|x_{s}}\,ds$. Following the averaging approach of \cite{Thompson2014} based on the separation of time scales of $x_{t}$ and $z_{t}$, $v_{t}$, we obtain the asymptotic behaviour for $0 < \epsilon \ll 1$ of the characteristic function for the joint system of $(v_{t}, z_{t})$ system, conditioned on a fixed value of $x_{t}$. The asymptotic solution demonstrates that the time-integral of $v_{t}$ converges in distribution to the integral of an $\alpha$-stable process, and we use this result to weakly approximate $x_{t}$ in (\ref{eqn:x_dynamics_stable_averaging}) by $X_{t}$, known as 
The process $U(t)$ that weakly approximates the slow dynamics of $u(t)$ satisfies   
\begin{equation}
dU_{t} = \bar{f}(U_{t})\,dt + \sigma \left( \frac{f_{2}(U_{t})}{g_{2}(U_{t})} \right)\diamond dL_{t}^{(\alpha,\beta)}, \label{eqn:N+approx_general}
\end{equation}
where $\bar{f}(U) = f_{1}(U) + [f_{2}(U)g_{1}(U)/g_{2}(U)]$.  The `$\diamond$' symbol indicates the Marcus interpretation of the stochastic differential terms \cite{Marcus1978}, which is analogous to the Stratonovich interpretation for Gaussian white noise driven systems.  If $\sigma f_{2}/g_{2}$ is a constant, then (\ref{eqn:N+approx_general}) reduces to the It\={o} interpretation. If $\alpha = 2$, then (\ref{eqn:N+approx_general}) is interpreted in the sense of Stratonovich. Details of theory and simulations are given in the 
Section \ref{nonlin2} and the references \cite{Thompson2014, Chechkin2014, Cont2004}.

Applying the (N+) stochastic averaging approximation (\ref{eqn:N+approx_general}) to the $(\eta,z)$ system defined in (\ref{eqn:dynamics_eta_with_z}) and (\ref{eqn:zt_dynamics}),  we obtain the SDE for the weak approximation $\tilde{\eta}_{t}$ 
\begin{equation}
d\tilde{\eta}_{t} = \tilde{f}(\tilde{\eta}_{t})\,dt + \epsilon^{\gamma^{*} - \rho}\Sigma dL_{t}^{(\alpha^{*},\beta^{*})}.
\label{eqn:reduced_linear_CAM_system_eta}
\end{equation}
Taking the inverse of the transformation $\mathcal{U}$ as in \cite{Thompson2014}, we obtain the SDE for $X_{t} = \mathcal{U}^{-1}(\tilde{\eta}_{t})$  that weakly approximates $x_{t}$ 
\begin{equation}
dX_{t} = f_{1}(X_{t})\,dt + \epsilon^{\gamma^{*} - \rho}\Sigma f_{2}(X_{t})\diamond  dL_{t}^{(\alpha^{*},\beta^{*})}. \label{eqn:reduced_linear_CAM_system}
\end{equation}
The value for $\Sigma$ used in this approximation is determined from  the relationship (\ref{eqn:int_y_to_stable}) and
is approximated numerically as discussed in
Appendix \ref{sec:estim_sigma}. 

\section{Sample systems}
\label{sec:sample_systems}
In this section, we simulate one linear and two nonlinear systems of the form (\ref{eqn:slowSDE}), (\ref{eqn:CAMsde}) and compare the stationary PDFs and ACDs of the numerically simulated trajectories of $x_{t}$ to those of the corresponding stochastic averaging approximation, $X_{t}$ (\ref{eqn:reduced_linear_CAM_system}). We take $\rho = \gamma^{*}$ for convenience, in which case the
power of $\epsilon$ in the noise coefficient of (\ref{eqn:reduced_linear_CAM_system}) vanishes.
%The details for numerically estimating $\Sigma$ are given in Appendix \ref{sec:estim_sigma}.

% We note here that when estimating the value of $\theta$ used in our simulations of $X_{t}$, we use $\epsilon = 10^{-5}$ which is less than $\epsilon$ considered in any of the full systems by at least one order of magnitude. Such a small value of $\epsilon$ is chosen to ensure that, asymptotically, the estimates of $\theta$ have no dependence (or at least a lessened dependence) on $\epsilon$ so as to emphasize the differences between the full and reduced systems due to the large, but finite, time scale separation.

\subsection{Linear system}
\begin{figure*}[t]
\centering
\includegraphics[width=0.45\linewidth]{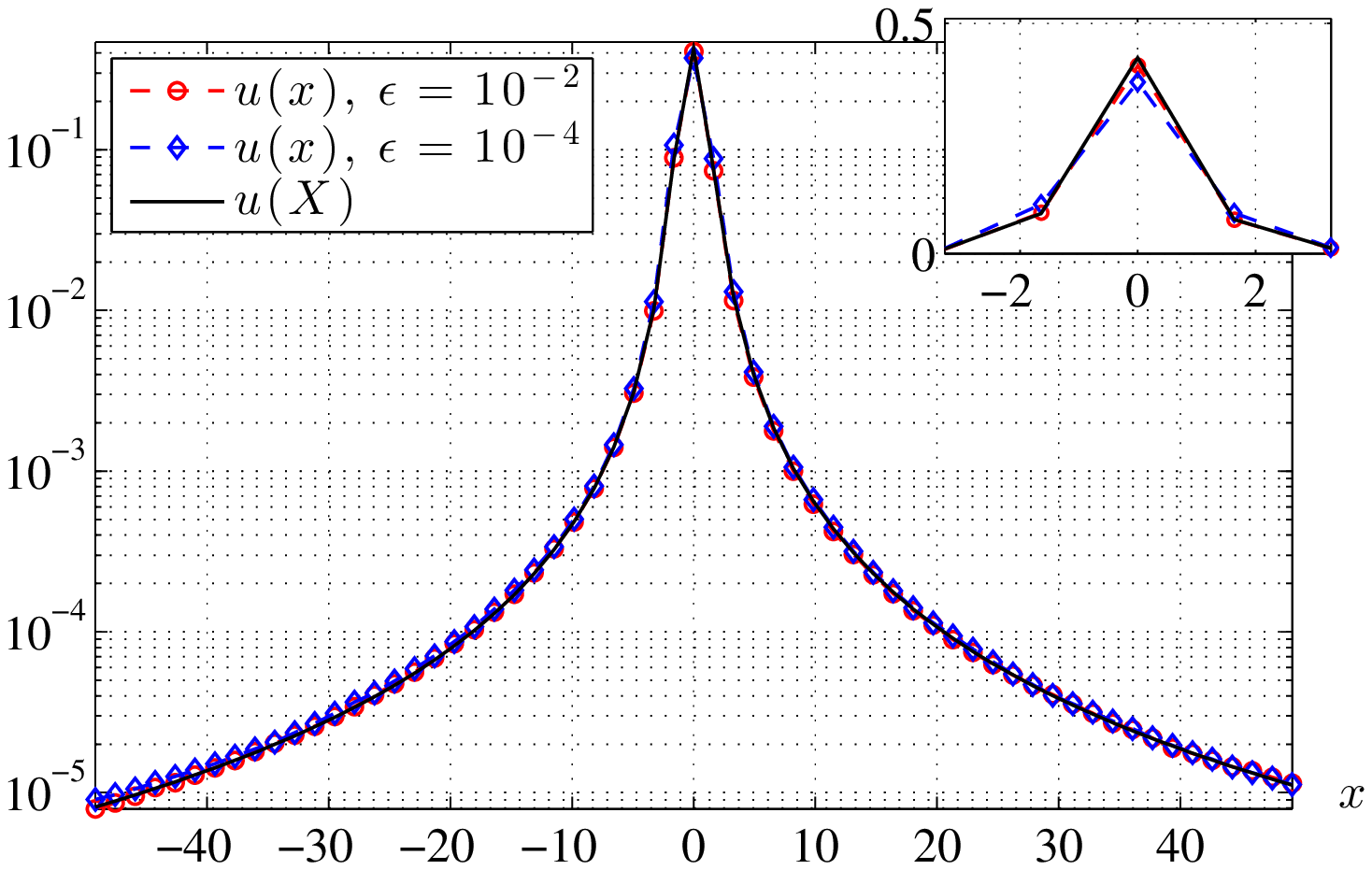}
\includegraphics[width=0.45\linewidth]{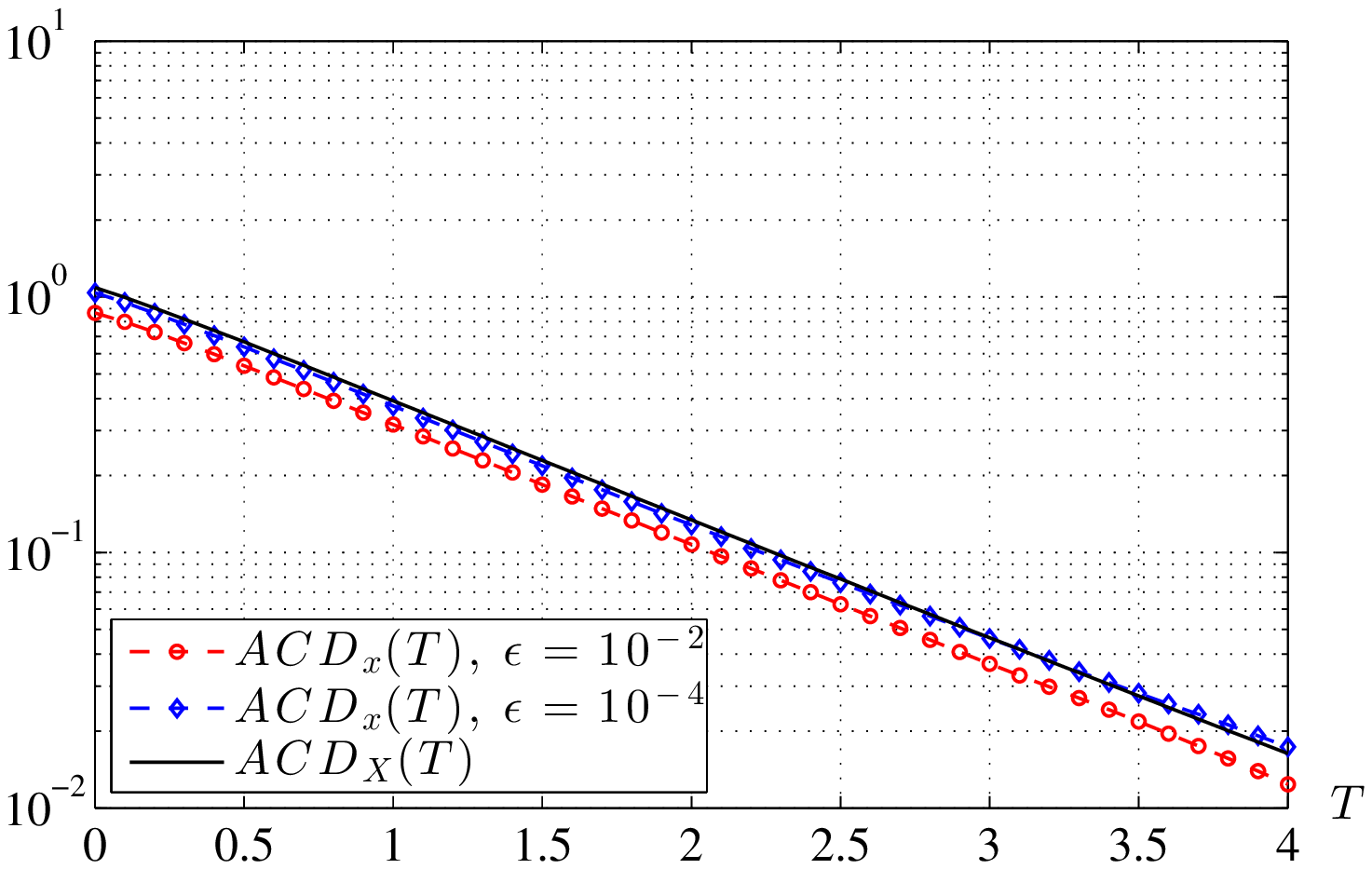} 
\includegraphics[width=0.45\linewidth]{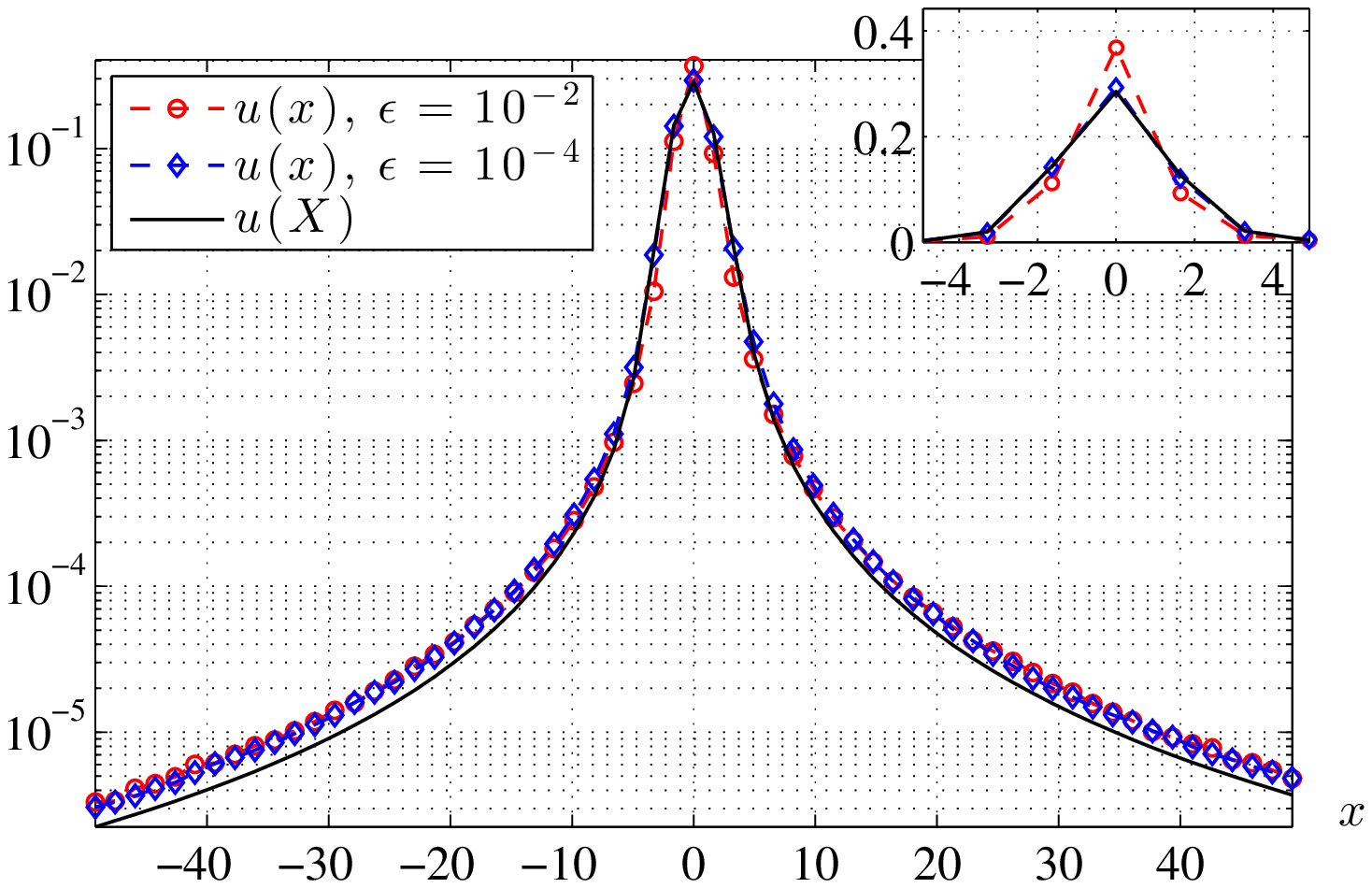}
\includegraphics[width=0.45\linewidth]{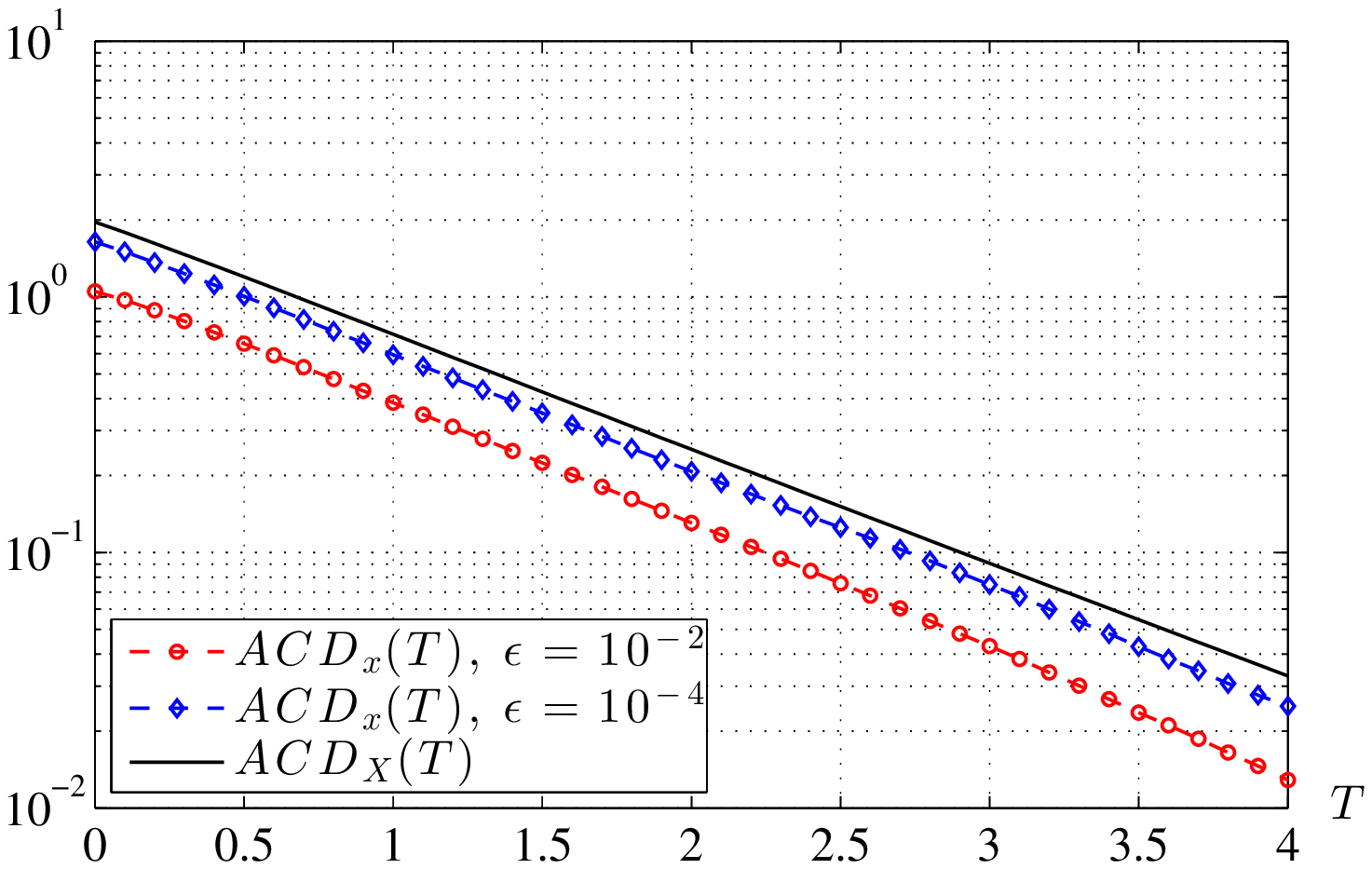}
\caption{{Top-Left}: The numerically estimated stationary PDF of $x_{t}$ (\ref{eqn:dynamics_x}) and the predicted PDF of $X_{t}$ (\ref{eqn:reduced_dynamics_x}) on a logarithmic scale (linear scale, inset) with parameters $(L,\alpha^{*},g,b) = (-1,1.5,0.1,0.5)$, $(\mu,\zeta) = (1,1)$. The value of $\epsilon$ is indicated in the legend. {Top-Right}: The corresponding ACD. {Bottom}: Same as top, but with $\alpha^{*} = 1.8$.  Note that the convergence rate of the reduced model toward the full system is smaller for the value of $\alpha^{*}$ closer to 2.}
\label{fig:lin}
\end{figure*}

\begin{figure}[t]
\centering
\includegraphics[width=0.49\linewidth]{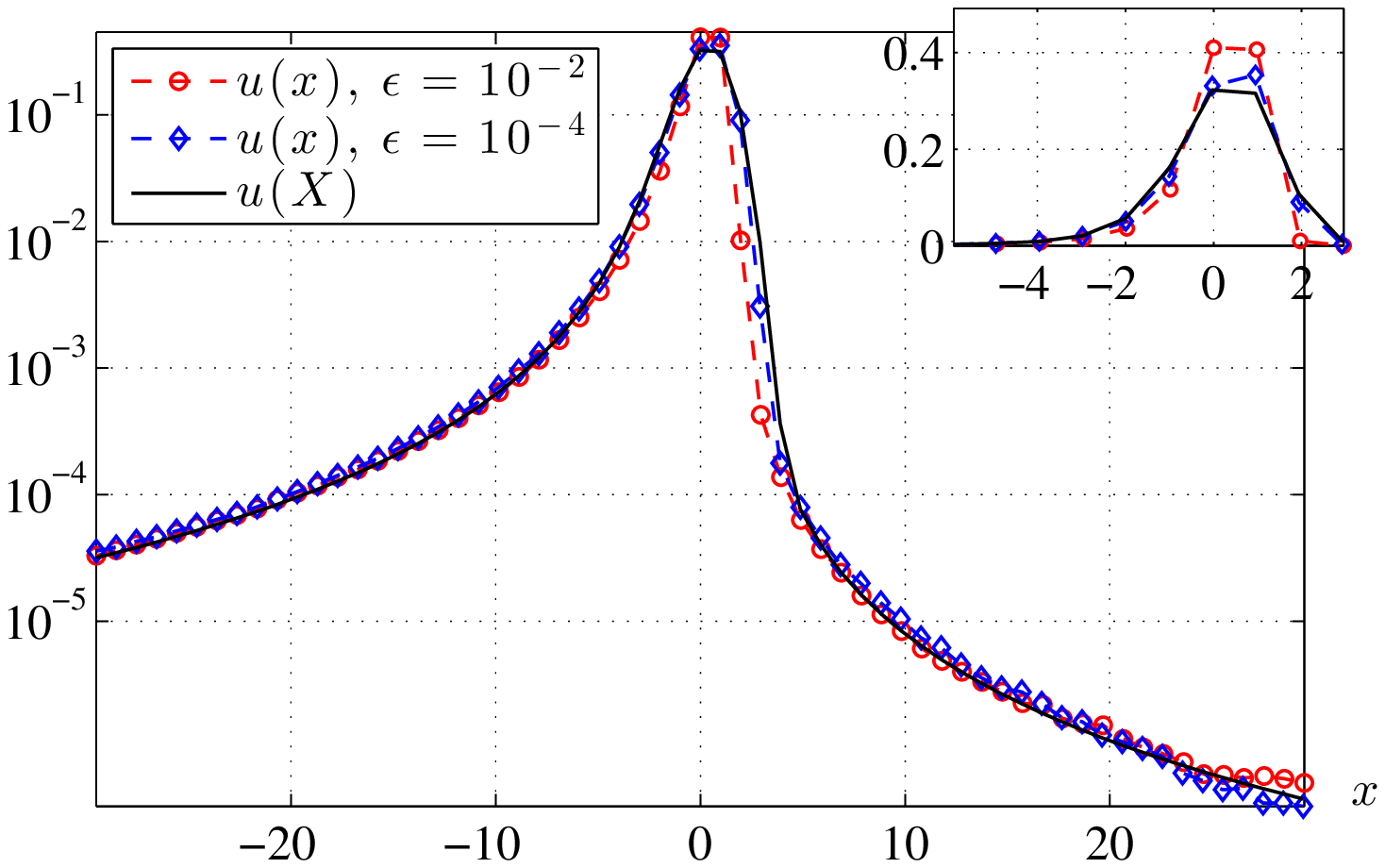}
\includegraphics[width=0.49\linewidth]{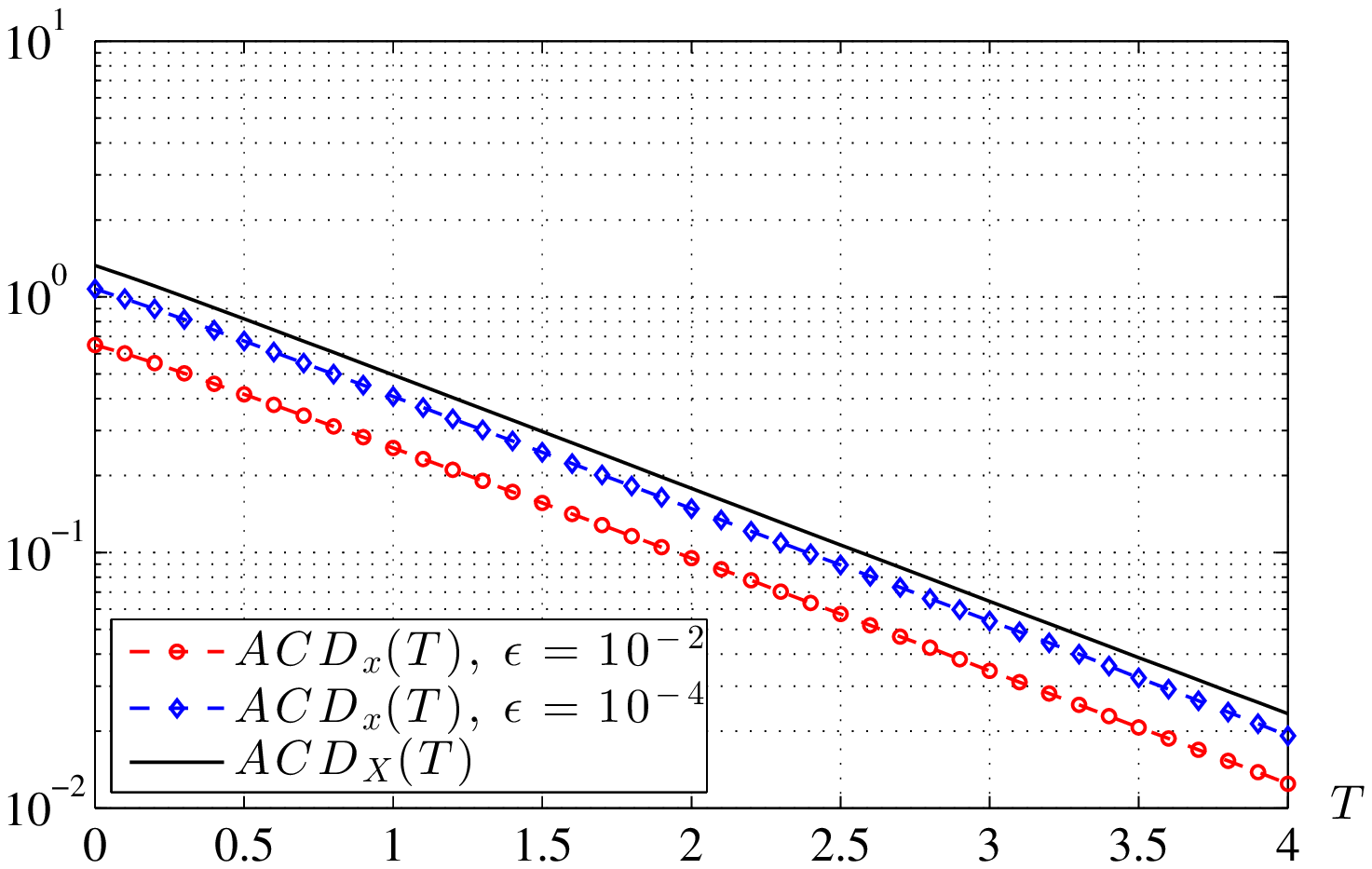}
\caption{As in Figure \ref{fig:lin}, but with parameters $(L,\alpha^{*},g,b) = (-1,1.7,-0.5,0.25)$, $(\mu,\zeta) = (1,1)$.}
\label{fig:lin2}
\end{figure}

Let the dynamics of $x_{t},\,t \ge 0$ be given by
\begin{equation}
dx_{t} = \left(-\mu x_{t} +  \epsilon^{-\gamma^{*}}\zeta y_{t/\epsilon}\right)\,dt \label{eqn:dynamics_x}
\end{equation}
where $\mu > 0$, $\zeta$ is a constant, $x_{0}$ is known and $y_{t/\epsilon}$ is a fast linear CAM noise process with dynamics (\ref{eqn:CAMsde}). Our stochastic averaging result (\ref{eqn:reduced_linear_CAM_system}) gives the reduced system,
\begin{equation}
dX_{t} = -\mu X_{t}\,dt +\zeta\Sigma  dL_{t}^{(\alpha^{*},\beta^{*})}. \label{eqn:reduced_dynamics_x}
\end{equation}
where $\sigma^{*}$ is given in (\ref{eqn:sigma_conv}).  Not surprisingly, given the linear form for (\ref{eqn:dynamics_x}), $X_{t}$ is an OULP.  
We simulate $x_{t}$ numerically and compare its estimated PDF and ACD functions with the known PDF and ACD for the corresponding process $X_{t}$. The results are shown in Figure \ref{fig:lin}   for  small $g$ (and thus small skewness parameter
$\beta^*$) and in Figure \ref{fig:lin2} for larger $g \ne 0$. It is evident in both cases that the stationary distribution and ACD of the reduced model converge to that of the full system as $\epsilon$ decreases.
% Notice that in Figure \ref{fig:lin}, the ACD functions of $x_{t}$ and $X_{t}$ do not match exactly for the $\alpha^{*} = 1.8$ case, but they appear to be the same up to a multiplicative constant. {\bf  I think Will is redoing these figures normalizing according to Adam's suggestion, so this discussion may not be relevant? This difference can be explained by considering the slow rate of convergence to a stable distribution for $\alpha^{*}$ near 2 and the fact that our estimator for $\Sigma$ assumes that the random variables $\int_0^{\Delta} y_{t/\epsilon} dt =Y_{j}$ have converged in distribution to an $\alpha$-stable distribution, which may not be the case. }{\bf Why would this shift the value of the ACD to a different value with same slope? Nevertheless, the differences between the slow variable of the full system and the reduced system are clearly decreasing as $\epsilon$ is reduced.  }

\subsection{Nonlinear system 1: nonlinear potential}
\begin{figure}[t]
\centering
\includegraphics[width=0.49\linewidth]{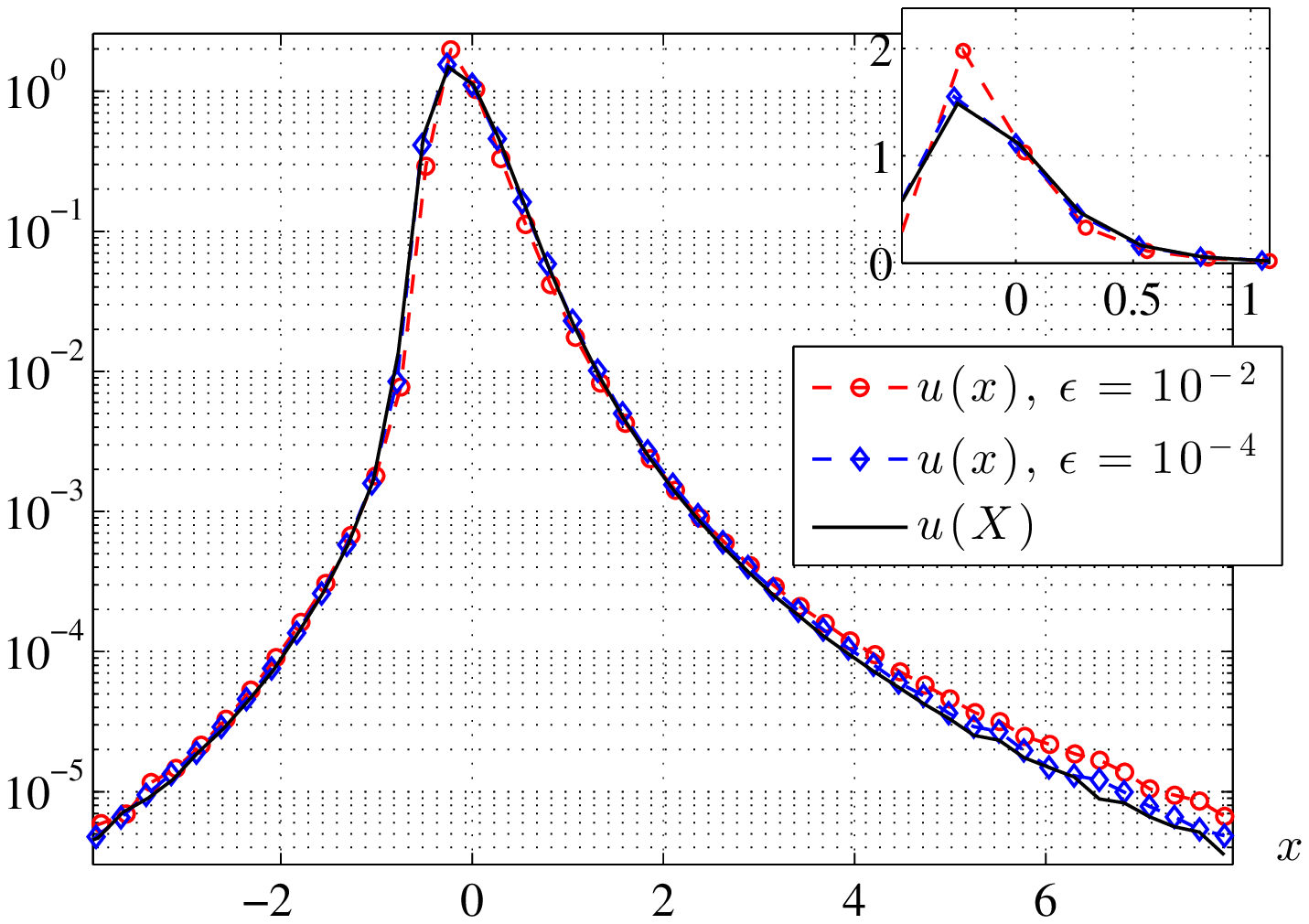}
\includegraphics[width=0.49\linewidth]{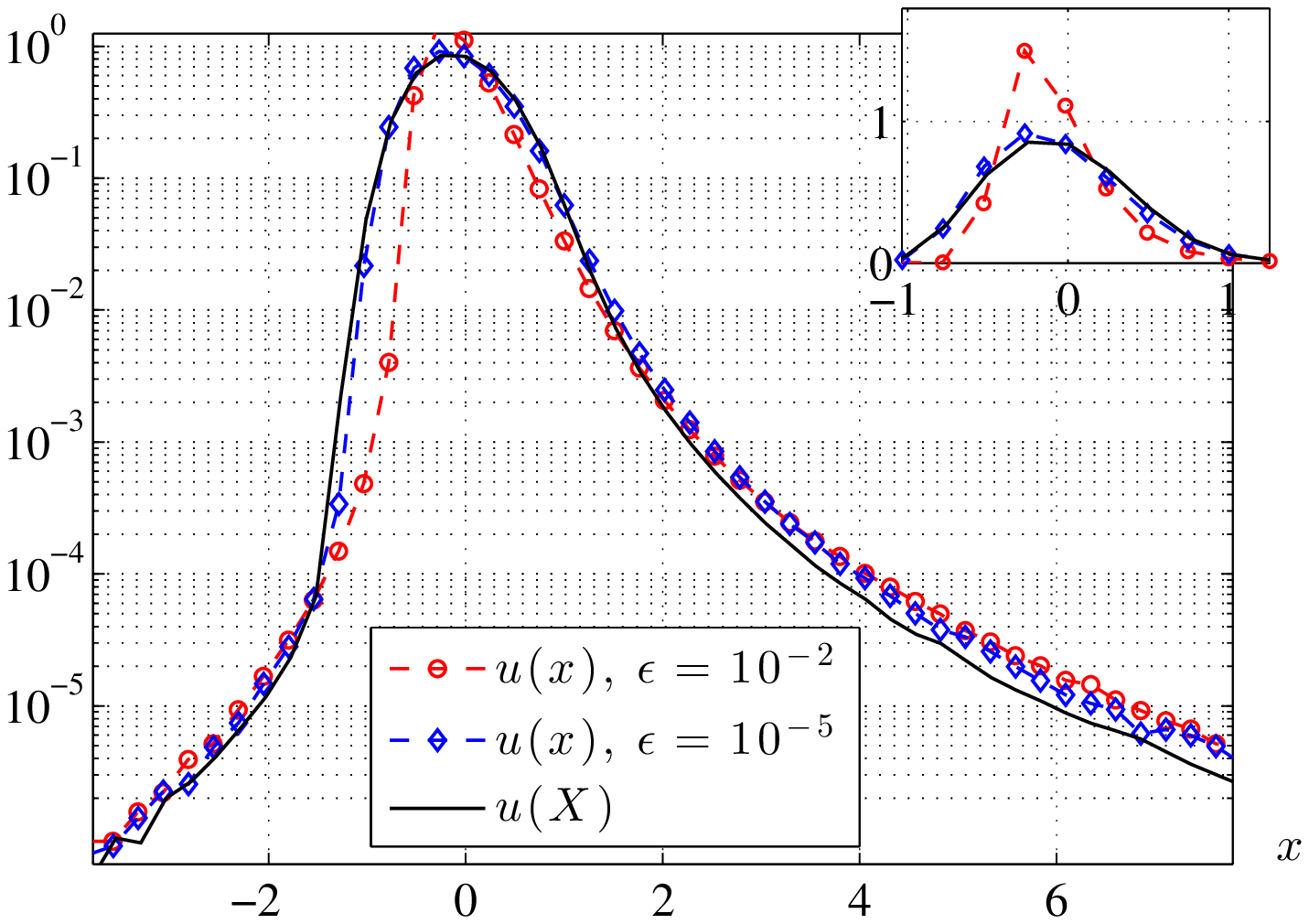}
\caption{The numerically estimated stationary PDFs $u(x)$ and $u(X)$ for nonlinear system 1 with dynamics (\ref{eqn:nonlin1}) and (\ref{eqn:nonlin1_reduced}), respectively. Results are shown on a logarithmic scale and a linear scale (inset). $(L,g,b) = (-1,1,0.5)$, $(\mu,\zeta) = (1,0.2)$ and $\epsilon$ as indicated in the legend. {Left/Right}: $\alpha^{*} = 1.5/1.8$.}
\label{fig:nonlin1}
\end{figure}

\begin{figure}[t]
\centering
\includegraphics[width=0.49\linewidth]{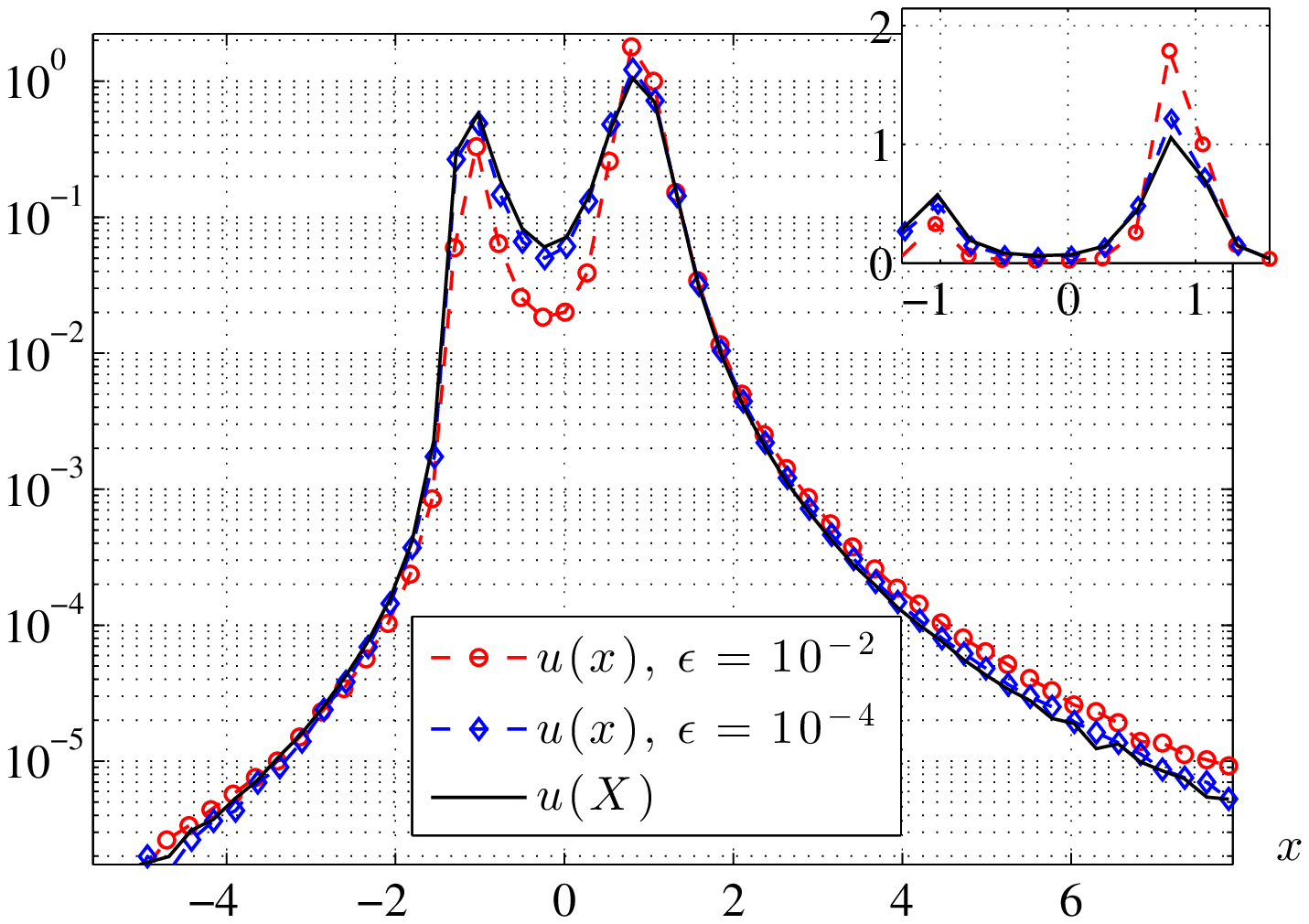}
\includegraphics[width=0.49\linewidth]{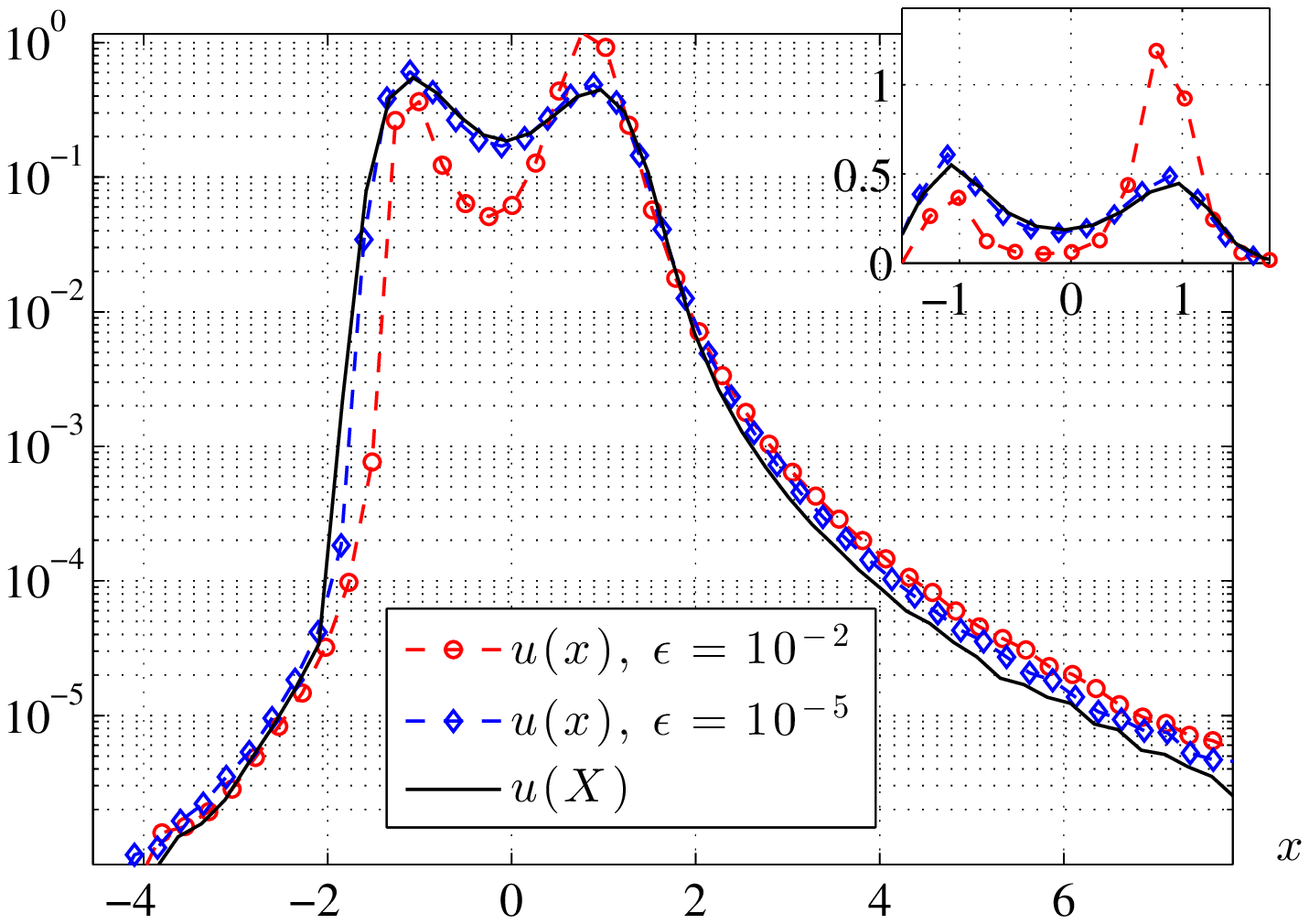}
\caption{  As in Figure \ref{fig:nonlin1}, but with parameters $(L, g,b) = (-1,1,0.5)$, $(\mu,\zeta) = (-1,0.2)$ and $\epsilon$ as indicated in the legend. {Left/Right}: $\alpha^{*} = 1.5/1.8$.}
\label{fig:nonlin1_multimode}
\end{figure}

We illustrate our reduction method in a system with cubic nonlinearity in the slow equation and additive linear CAM noise driving:
\begin{equation}
dx_{t} = -\left(\mu x_{t} + x_{t}^{3}\right)\,dt + \epsilon^{-\gamma^{*}}\zeta\,y_{t/\epsilon}\,dt \label{eqn:nonlin1}
\end{equation}
where $y_{t}$ is given by (\ref{eqn:CAMsde}). This system is close to linear for small $x_{t}$, but experiences a stronger nonlinear drift for large $x_{t}$. For $\mu > 0$, the system is globally attracted to $x = 0$ in the noise-free limit, and for $\mu < 0$, the origin is an unstable (repelling) equilibrium and the system is locally attracted to one of two stable equilibria at $x = \pm \sqrt{-\mu}$. According to our reduction results, we expect that the dynamics for $x_{t}$ can be weakly approximated by
\begin{equation}
dX_{t} = -(\mu X_{t} + X_{t}^{3})\,dt + \zeta \Sigma\,dL_{t}^{(\alpha^{*},\beta*)}, \label{eqn:nonlin1_reduced}
\end{equation}
where $\Sigma$ is numerically estimated in the same fashion as described in Section \ref{subsec:CAM2stablenoise}. We use a predictor-corrector method described in \cite{Thompson2014} to simulate both the full and reduced systems, avoiding numerical instabilities due to the cubic nonlinearity for large values of $x_t$ or $X_t$.   The CAM noise process is integrated as described in Appendix \ref{app:simulating_CAM}.  The numerically estimated PDFs for $\mu > 0$ and $\mu<0$ are compared in Figures \ref{fig:nonlin1} and \ref{fig:nonlin1_multimode}, respectively.  In both cases, the stochastic averaging approximation works well for sufficiently small $\epsilon$. 
The PDF tails from both the full and reduced systems show modest fluctuations due to sampling variability.

\subsection{Nonlinear system 2: bilinear interaction}\label{nonlin2}
\begin{figure}[t]
\centering
\includegraphics[width=0.49\linewidth]{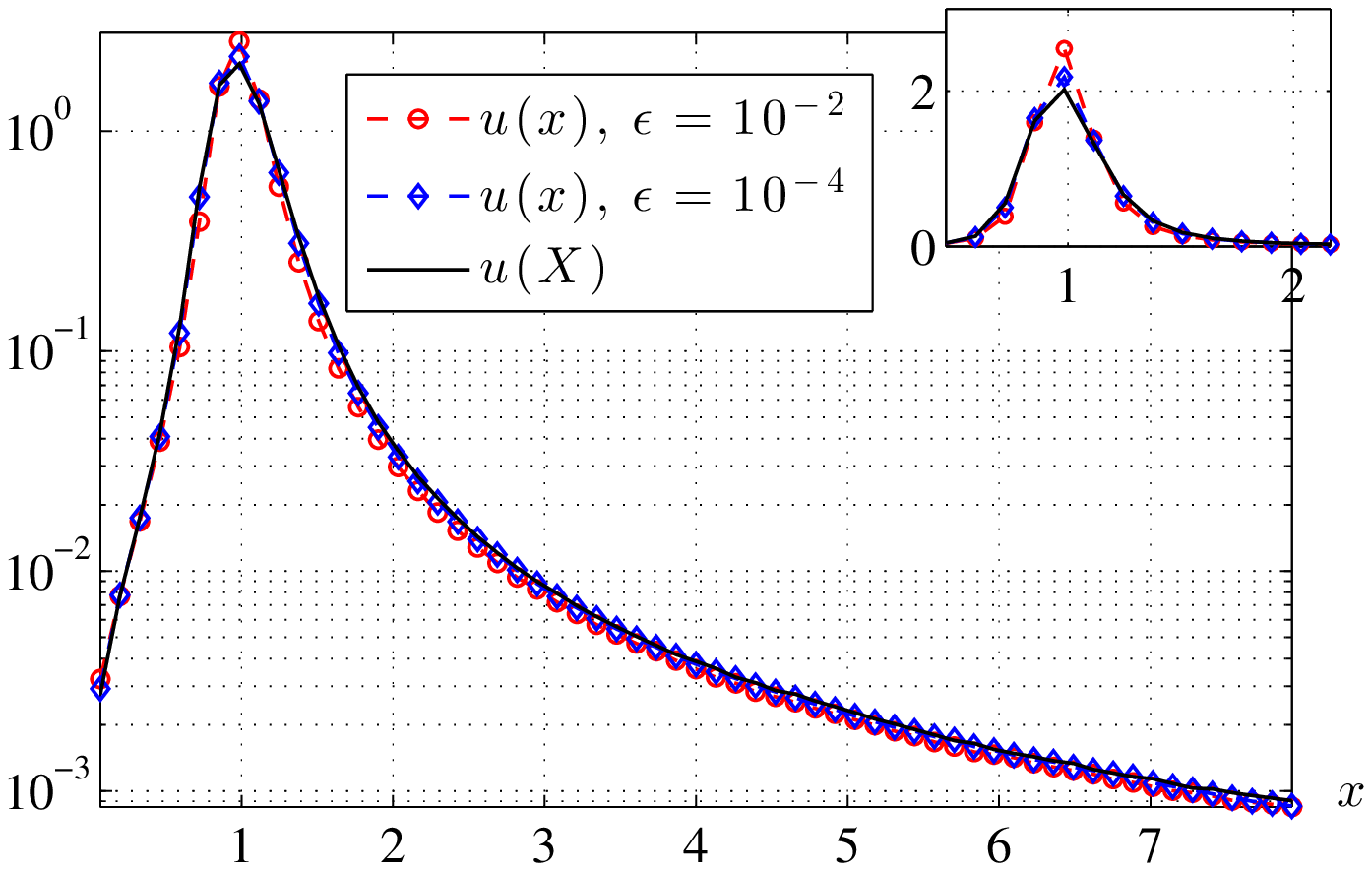}
\includegraphics[width=0.49\linewidth]{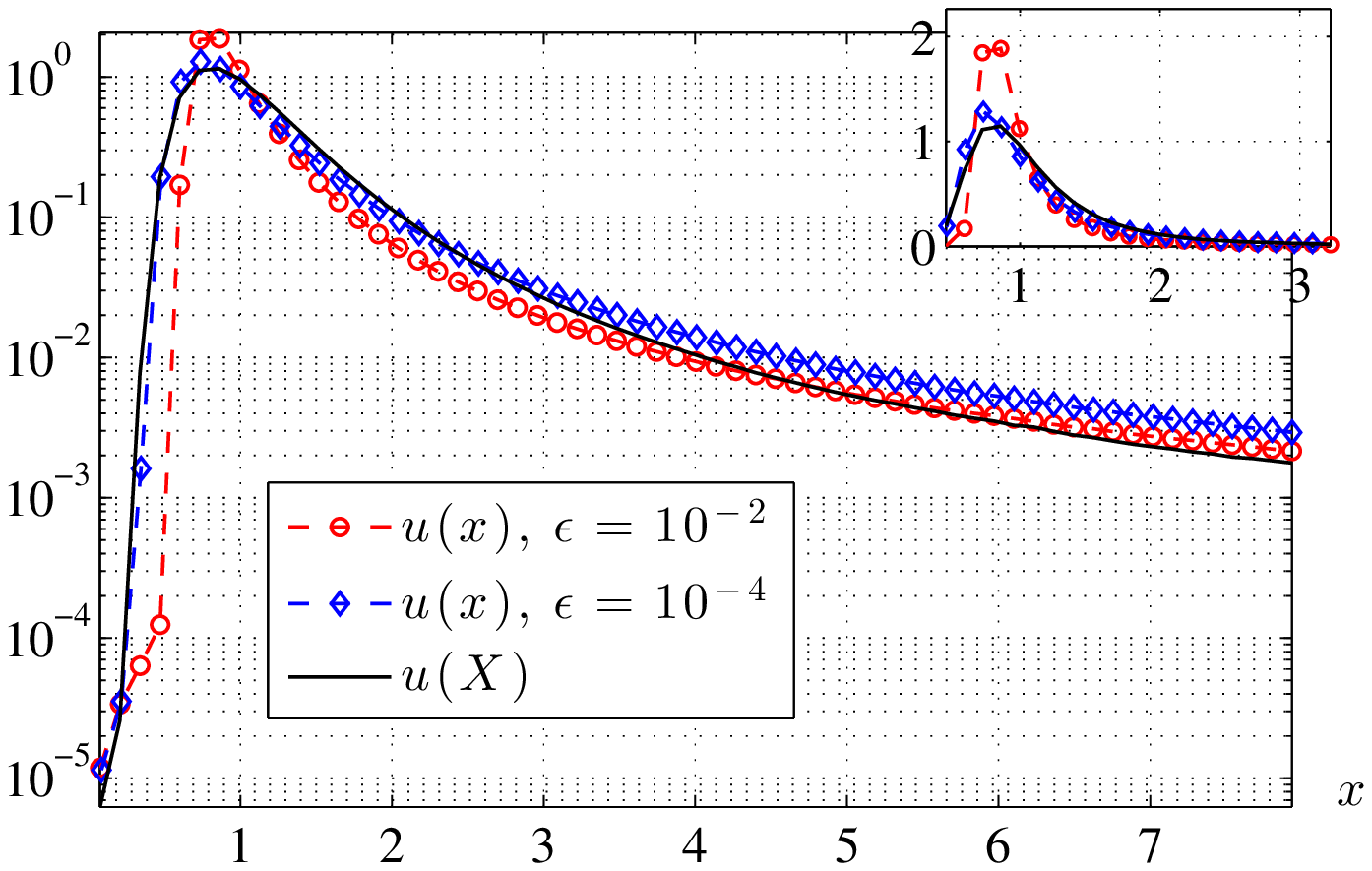}
\caption{The numerically estimated stationary PDFs $u(x)$ and $u(X)$ for nonlinear system 2 with dynamics (\ref{eqn:nonlin2}) and (\ref{eqn:nonlin2_reduced}), respectively, shown on a logarighmic scale and linear scale (inset).. Left:  $(L,\alpha, g,b) = (-1,1.5,0.1,0.5)$, $(c,\zeta) = (1,0.2)$. Right: $(L,\alpha, g,b) = (-1,1.8,1,0.1)$, $(c,\zeta) = (1,0.2)$.}
\label{fig:nonlin2}
\end{figure}

The second nonlinear system that we consider has a bilinear term in which the CAM noise process appears multiplicatively in the slow
equation:
\begin{equation}
dx_{t} = (c - x_{t} + \zeta x_{t}y_{t/\epsilon})\,dt, \quad c,\,x_{0} > 0. \label{eqn:nonlin2}
\end{equation}
It is clear from considering the limit as $x \to 0$ that the dynamics of $x_{t}$ are restricted to the positive real line, due to the inclusion of the state-independent drift term $c > 0$. Applying our stochastic averaging result, we can weakly approximate the dynamics of (\ref{eqn:nonlin2}) by $X_{t}$ where
\begin{equation}
dX_{t} = (c - X_{t} )\,dt + \zeta \Sigma\, X_{t} \diamond dL_{t}^{(\alpha^{*},\beta^*)}, \quad X_{0} = x_{0}.\label{eqn:nonlin2_reduced}
\end{equation}
where `$\diamond$' denotes the Marcus interpretation of the stochastic driving term. Details and references on the calculation of Marcus stochastic integrals
used in simulating these SDE's can be found in \cite{Thompson2014, Chechkin2014, Cont2004}. To summarize,  the Marcus stochastic term, $\zeta \Sigma X_{t}\diamond dL_{t}^{(\alpha^{*},\beta^*)}$,  can be written as $\lambda(1;\Delta L, X_{t-}) - X_{t-}$ where $\lambda(s;\Delta L, X_{t-}) = \lambda(s)$ satisfies
\begin{equation}
\frac{d\lambda(s)}{ds} = \zeta \Sigma(\Delta L)\, \lambda(s), \quad \lambda(0) = X_{t-}, \label{eqn:nonlinsys2_Marcus_incre}
\end{equation}
with $\Delta L = dL_{t}^{(\alpha^{*},\beta^*)}$ is the jump of the $\alpha$-stable process at time $t$ and $X_{t-} = \lim_{s \to t^{-}}X_{s}$. 
In this case, (\ref{eqn:nonlinsys2_Marcus_incre}) can be solved for $\lambda(1)$, giving $\lambda(1;\Delta L, X_{t-}) = \exp(\zeta \Sigma\, \Delta L)X_{t-}$. Having an analytic expression for the Marcus stochastic increment facilitates a straightforward and accurate simulation of (\ref{eqn:nonlin2_reduced}). The estimates of the stationary PDFs of $x_{t}$ and  the reduced process $X_{t}$ compare well (Figure \ref{fig:nonlin2}).

\subsection{Discussion of numerical results}\label{discuss_num}
From the results displayed in Figures \ref{fig:lin} - \ref{fig:nonlin2}, we see that the stationary behaviour of $x_{t}$ in each system under consideration can be well-approximated by a system influenced by $\alpha$-stable noise for sufficiently small values of $\epsilon$. Both the PDFs and the ACD functions of the full and reduced systems match well. 

As expected for stochastic averaging approximations for systems with multiple time scales, as $\epsilon$ decreases, we see improvement of the approximation of the weak properties of $x_t$ in  the full system by $X(t)$ in the reduced system. However, the difference between the PDFs and ACDs of the reduced systems and the full ones depends not only on $\epsilon$, but also on $\alpha^{*}$. Specifically,  for $\alpha^{*}$ close to 2 the value of $\epsilon$ needs to be noticeably smaller, as we can see by reviewing Conditions A) and B) from Section \ref{subsec:CAM2stablenoise}  for  $\int_0^T y_{s/\epsilon} ds$ to be approximately distributed as $\alpha$-stable.  A sum of independent heavy-tailed random variables converges to an $\alpha$-stable variable more slowly for values of $\alpha$ closer to 2 \cite{Kuske2001}.  As such, the number of approximately independent terms in the discretization of the integral $\int_{0}^{T}y_{s/\epsilon}ds$ for fixed $T$ must be larger for $\alpha$ closer to two in order for this integral to have the desired properties - requiring a smaller value of $\epsilon$.   The quality of the weak approximation of $x_{t}$ by $X_{t}$ is determined not only by the timescale separation, but also by the tail behaviour of the limiting $\alpha$-stable process.

\section{Conclusions}
\label{sec:conclusions}

This study demonstrates that $\alpha$-stable forcing can appear in the asymptotic slow dynamics of fast-slow systems where the fast process is linear and forced with a combination of additive and multiplicative noise such that its distribution has infinite variance. We studied the linear CAM process in this analysis due to its use in applied research problems and the fact that many analytical results are available \cite{Sardeshmukh2009, Penland2012}. A particularly valuable aspect of the linear CAM process is that it provides a simple dynamical form in which the interaction of multiplicative terms and Gaussian white noise forcing results in a process possessing a distribution with power law tails and infinite variance. Due to the infinite variance stationary behaviour of the linear CAM noise process and the explicit predictions of the power law tail behaviour, we hypothesized that the linear CAM noise process would appear to be equivalent to an $\alpha$-stable forcing term when used to drive a slower process. The derivation of the approximation required us to focus our attention on the properties of the integral of the linear CAM noise process and determining an OULP  having integral behaviour with similar statistics. Formulas for the stability index and skewness parameter of the corresponding OULP in terms of the CAM process parameters are straightforward. In contrast, a coefficient in the OULP drift term needs to be estimated numerically as it depends on the serial dependence properties of the linear CAM noise process for which analytic results are not available. When we apply our stochastic averaging approximation to linear systems as well as systems with nonlinearities in the slow variable, we observe good agreement between simulations of the full and reduced dynamics, which improve as $\epsilon$ decreases. For values of the parameters such that the 
%tails of the stationary distribution of the linear CAM process has decay exponent close to 3, 
corresponding OULP has a stability index close to 2, 
the ratio of time scales between the fast and slow processes needs to be extremely large to observe the distributional convergence of the slow variable to the predicted distribution.

Besides presenting a method of approximating fast-slow systems that are forced with correlated additive and multiplicative noise processes, this analysis also suggests a possible mechanism through which $\alpha$-stable forcing can emerge in the modelling of physical problems. For example, as mentioned above, CAM noise processes emerge from considerations of the dynamics of quadratically nonlinear systems like those describing atmospheric motion \cite{Majda2001, Sardeshmukh2009} and could result on longer time scale processes experiencing forcing terms distributed according to an approximately $\alpha$-stable law. This result offers a possible explanation for the observation of $\alpha$-stable noise in %the calcium signal found in %the Greenland ice cores \cite{Ditlevsen1999}.
 various climatic and  fluid dynamical time series \cite{Huber2001,Seo2000,Ditlevsen1999,Del-Castillo-Negrete1998}.

% FUTURE WORK
	% Multi-dimensional systems
	% More general nonlinearities
	% Other infinite variance processes
	% Applications to real world systems
There are various extensions to this research that are worth exploring. All of the analysis given is for the case where both the fast and slow subsystems are univariate. Extensions of the results of this paper to higher-dimensional fast slow systems may yield unexpected challenges, but are required before this stochastic averaging method could be applied to more general research problems. 
%To have explicit results showing that a climate model can exhibit $\alpha$-stable statistics would certainly be of interest and \cite{Huber2001, Seo2000, Ditlevsen1999,Del-Castillo-Negrete1998} provide some examples for the appearance of infinite variance statistics in climate processes. 
Also, while there is broad utility in studying models where the fast linear CAM noise process perturbs the slow variable, it would be worth exploring the implications of nonlinear fast perturbations to the slow process or where the fast, infinite-variance process is something other than the linear CAM noise process. Another interesting situation would be to study the situation where the fast linear CAM process is conditionally dependent of the slow process, or in other words the slow variable influences the parameters of the linear CAM noise process. This would be a particularly useful scenario to consider in the context of climate modelling.

 \begin{acknowledgements}
AHM and RK acknowledge partial support from the Natural Sciences and Engineering Research Council (NSERC) Discovery Grant program. RK
also was partially supported by a grant from the Simons Foundation for work carried out in part at the Isaac Newton Institute.  WFT acknowledges support from the NSERC Alexander Graham Bell Canadian Graduate Scholarship program and the UBC Faculty of Science Graduate Award program.
\end{acknowledgements}

\appendix

\section{Characteristic functions}
\subsection{Calculations for $\psi$}\label{calc_details}

We provide some details of the calculations for the characteristic functions arising in Section \ref{subsec:OULP_approx_for_y}.

Taking the Fourier transform of (\ref{eqn:FKE_vz}), gives a quasi-linear partial differential equation whose solution is the joint characteristic function  $\psi_{v,z}(m,k,t) = \mathcal{F}[P](m,k,t) = \iint_{\mathbb{R}^{2}} \exp(ikz + imv) P(v,z,t)\,dvdz$.
\begin{eqnarray}
\frac{\dd \psi_{v,z}}{\dd t} + \left(\frac{\theta k}{\epsilon} - m\right)\frac{\dd \psi_{v,z}}{\dd k} &=& -\frac{\left(\sigma_z\right)^{\alpha}}{\epsilon}|k|^{\alpha^{*}}\Xi(k;\alpha,\beta)\psi_{v,z}, \label{eqn:FT_FKE_vz} \\
\psi_{v,z}(m,k,0) &=& \exp(ikz_{0}) 
\end{eqnarray}
where $\Xi(\Lambda;\alpha,\beta) = 1 - i\beta\sgn{\Lambda}\tan(\pi\alpha/2)$ as given in (\ref{eqn:Xi}). Solving (\ref{eqn:FT_FKE_vz}) via the approach  in \cite{Thompson2014}, the method of characteristics gives the solution for $\psi_{v,z}$
\begin{equation}
 \psi_{v,z}(m,k,t) = \exp\left(ikz_{0}e^{-\theta t/\epsilon} + i\frac{\epsilon m z_{0}}{\theta}(1 - e^{-\theta t/\epsilon})  - \frac{(\sigma_z)^{\alpha}}{\epsilon}\int_{0}^{t}|\Lambda(r)|^{\alpha}\Xi(\Lambda(r);\alpha,\beta)\,dr\right) \label{eqn:char_func_joint_vz_sys}
\end{equation}
where $\Lambda(r) = \frac{\epsilon m}{\theta} + \left(k - \frac{\epsilon m}{\theta}\right)e^{-\theta r /\epsilon}$.  Following \cite{Thompson2014},  the
 integral term in (\ref{eqn:char_func_joint_vz_sys}) has the asymptotic behavior for $t = O(1)$ and $\epsilon\ll 1$,
\begin{align}
\frac{(\sigma_z)^{\alpha}}{\epsilon}\int_{0}^{t}|\Lambda(r)|^{\alpha}&\Xi(\Lambda(r);\alpha,\beta)\,dr \label{eqn:asymp_integral} \\ &= \left(\frac{\sigma_z^{\alpha}}{\alpha^{*}\theta}\right)|k|^{\alpha}\Xi(k;\alpha,\beta) + \frac{\epsilon^{\alpha-1}\sigma_z^{\alpha} t}{\theta^{\alpha}}|m|^{\alpha}\Xi(m;\alpha,\beta) + O(\epsilon). \nonumber
\end{align}
Note that this expression can be factored into separate functions of $k$ and $m$ for $\theta t /\epsilon\gg 1$, resulting in the expressions for
$\psi_v$ and $\psi_z$ in (\ref{eqn:char_func_z} - \ref{eqn:char_func_v}).

\subsection{Asymptotic behaviour of the characteristic function $\psi_Y(k)$ for small $k$}
\label{subsec:asymp_charfunc}

The characteristic function of $Y_j$ is given by the Fourier transform, $\psi_{Y}(k)$, of the PDF $u_{Y}(r)$ in (\ref{eqn:dist_RR}) 
\begin{align}
\psi_{Y}(k) &= \int_{\mathbb{R}}\exp(iky)u_{Y}(y)\,dy = q^{-}\left(\int_{-\infty}^{-a} \frac{e^{iky}\,dy}{|y|^{-(1 + \alpha^*)}}\right) + \left(\int_{-a}^{a}e^{iky}u_{Y}(y)\,dy\right) + q^{+}\left(\int_{a}^{\infty} \frac{e^{iky}\,dy}{y^{-(1 + \alpha^*)}}\right) \\
&= \left( \Gamma(-\alpha) (q^{+} + q^{-})\cos\left(\frac{\pi\alpha^*}{2}\right)\right)|k|^{\alpha^*}\Xi(k;\alpha^*,{\beta_Y})
%\left(1 - i\tilde{\beta}\sgn{k}\tan\left(\frac{\pi\alpha}{2}\right)\right)
   \label{eqn:intermed1_1} \\  & \quad + q^{-}\left(\int_{a}^{\infty}\frac{1 - iky}{y^{1 + \alpha^*}}\,dy - \sum_{n = 2}^{\infty}\frac{(-ika)^{n}a^{-\alpha^*}}{(n - \alpha^*)\Gamma(n+1)} \right) + q^{+}\left(\int_{a}^{\infty}\frac{1 + iky}{y^{1 + \alpha^*}}\,dy -  \sum_{n = 2}^{\infty}\frac{(ika)^{n}a^{-\alpha^*}}{(n - \alpha^*)\Gamma(n+1)} \right) \nonumber \\ &\quad + \left(\int_{-a}^{a}e^{iky}u_{Y}(y)\,dy\right). \nonumber
\end{align}
Here $\Xi$ is given in (\ref{eqn:Xi}), $\beta_Y = \frac{q^{+} - q^{-}}{q^{+} + q^{-}}$, and we have used the facts $\int u_{Y}(r) dr = 1$ and $E[Y_{j}]=0$ to rewrite the integrals as in \cite{Kuske2001} in order to facilitate an
expansion for small $k$. Using the result from \cite{Sato1999}  
\begin{equation}\int_{b}^{\infty}\frac{e^{iky}\,dy}{y^{1 + \alpha^*}} = |k|^{\alpha^*}\Gamma(-\alpha^*)e^{-i\sgn{k}\pi\alpha^*/2} + \int_{b}^{\infty}\frac{1 + iky}{y^{1 + \alpha^*}}\,dy - \sum_{n = 2}^{\infty}\frac{(ik)^{n}b^{n-\alpha^*}}{(n - \alpha^*)\Gamma(n+1)}. \end{equation}
in an expansion for small $k$ as in \cite{Kuske2001}) yields
\begin{align}
\psi_{Y}(k) & = \exp\left[\left( \Gamma(-\alpha^{*}) (q^{+} + q^{-})\cos\left(\frac{\pi\alpha^{*}}{2}\right)\right)|k|^{\alpha^*}\Xi(k;\alpha^*,\beta_Y)
%\left(1 - i\beta^{*}\sgn{k}\tan\left(\frac{\pi\alpha^{*}}{2}\right) \right)
 - Qk^{2} + O(k^{3})\right]. \label{eqn:apprx_charfuncY}\\
Q &= \int_{-a}^{a}y^{2}u_{Y}(y)\,dy - \left(\frac{(q^{+} + q^{-})}{2 - \alpha^*}\frac{a^{2 - \alpha^*}}{2}\right). \label{eqn:Qfactor}
\end{align}

\section{The probability density $u_Y$}\label{app:uY}
For each of the $y_m$ in $Y_j \approx  \sum_m  y_{m} h$, $m= 1,2, \ldots M$, we use the Euler-Maruyama approximation for
(\ref{eqn:CAMsde}) to express $y_m$ in terms of $y_1$ and the random variables $\xi_{k,m} \sim N(0, \sqrt{h})$, $k=1,2$ in the approximation. Following this straightforward but tedious calculation, $Y_j$ takes the form
 \begin{eqnarray}
 Y_j &=&  y_1\left[ {\cal{C}}_1(\Delta)  + h \sum_{m=1}^{M} \sum_{\ell<m} \xi_{1,\ell} +
  h\sum_{m=1}^{M} \sum_{\ell\neq k<m }[ {\cal K}_1\xi_{1,\ell}
 \xi_{1,k} + {\cal K}_2 \xi_{1,\ell}\xi_{2,k}]\right] \nonumber\\
 & &\qquad +  h \sum_{m=1}^{M}\sum_{\ell<m}  (g\xi_{1,\ell}+ b\xi_{2,\ell}) \nonumber\\
 & \equiv & {\cal{C}}_1(\Delta)y_1 + h {\cal C}_2 y_1 +  h{\cal C}_3  = {\cal{G}}(y_1)\label{Yjtoy1}
 \end{eqnarray}
 where ${\cal C}_1$ is linear in $\Delta$, ${\cal{K}}_i$ for $i= 1,2$ are $O(1)$ constants dependent on $E$.
   We use (\ref{Yjtoy1}) in the expression for the density of $Y_j$ obtained via conditioning on  $y_1$ and $\xi_i,m$
\begin{eqnarray}
u_Y(Y_j) &=& \int\int ...\int  P(Y_j | y_1, \xi_{k,m}, m= 1, \ldots M; k=1,2) p_s(y_1)\prod_{m=1}^{M}\rho_1(\xi_{1,m})\rho_2(\xi_{2,m}) d\xi_{1,m}d\xi_{2,m}dy_1 \nonumber\\
&=& \left\langle\delta(Y_j - {\cal{G}}(y_1))p_s(y_1)\right\rangle_{\xi_{k,m}} \qquad   k = 1,2; \ m=1\ldots M
\label{eqn:Y_jpdf}
 \end{eqnarray}
 assuming $y_1$ is taken from the stationary distribution for $y$, $\rho_i$ is the density for $\xi_{i,m}$, and $\langle \cdot \rangle_{\xi}$ indicates expected value with respect to $\xi$.
 Using that ${\cal G}$ is  linear in $y_1$, and that  ${\cal C}_i$ for $i=2,3$ are  sums of products of independent random variables
 $\xi_{i,m}$,
 %\sim N(0,\sqrt{h})$, 
 we conclude that $u_Y$ has the behavior given in (\ref{eqn:dist_RR}).
 
\section{Simulating the CAM noise process}
\label{app:simulating_CAM}
% HOW DO WE SIMULATE THE CAM NOISE PROCESS?
 We use the weak order 2.0 explicit method \cite{Kloeden1992} to simulate the CAM noise process
  in  (\ref{eqn:CAMsde}). This weak numerical approximation takes the form for  $\hat{y}_{n} = y_{n\delta{t}}$,
\begin{align}
\hat{y}_{n+1} = \hat{y}_{n}& + \frac{\tilde{L}}{2}\left(\hat{y}_{n} + \Upsilon \right)\delta{t} + \frac{1}{4}\left(E\Upsilon^{+} + 2E\hat{y}_{n} + E\Upsilon^{-} + 4g\right)\,\delta{W_{1,n}} \nonumber\\ & + \frac{1}{4}\left( E\Upsilon^{+} - E\Upsilon^{-}\right)\left(\frac{\delta{W_{1,n}}^{2} - \delta{t}}{\sqrt{\delta{t}}} \right)+ b\,\delta{W_{2,n}} \label{eqn:weaksimscheme}
\end{align}
where
\begin{equation}
\begin{cases}
\Upsilon = \hat{y}_{n} + \tilde{L}\hat{y}_{n}\delta{t} + (E\hat{y}_{n} + g)\delta{W_{1,n}} + b\,\delta{W_{2,n}}, \\
\Upsilon^{\pm} = \hat{y}_{n} + \tilde{L}\hat{y}_{n}\delta{t} \pm (E\hat{y}_{n} + g)\sqrt{\delta{t}} + b\,\delta{W_{2,n}}.
\end{cases}
\end{equation}
Here $\tilde{L} = L + E^{2}/2$ and $\delta{t}$ is the size of the discrete time step. The terms $\delta{W_{1,n}},\,\delta{W_{2,n}},\, n = 0,1,2, \dots$ are independent Gaussian random variables with mean 0 and variance $\delta{t}$.

\subsection{Consistency of simulations}
% DO THE SIMULATIONS COMPARE WELL WITH THEORY? (when applicable)
Since $y_{t/\epsilon}$ evolves on the fast time scale, $\delta t$ must be chosen smaller than $\epsilon$ in order to resolve the fast dynamics. We compared the stationary PDF $p_s$ for $y_{t/\epsilon}$ given by (\ref{eqn:CAMpdf}) to numerical approximations to $p_s$  based on simulations of (\ref{eqn:CAMsde}).  As expected, $\delta t$ must be an order of magnitude smaller than $\epsilon$ in order to obtain relative
errors of the numerical simulations to $p_s$ that are $O(10^{-2})$ or smaller. We found that the relative error of two cases, $\delta{t} = \epsilon/10$ and $\epsilon/100$ were of the same order of magnitude. 
Furthermore, we also considered whether the value of $\delta t$ has an effect on the estimates of $Y_{j} = \int_{(j - 1)\Delta}^{j \Delta}y_{s/\epsilon}\,ds$. We compared the behaviour of the density of $Y_j$ for decreasing values of $\epsilon$. As in the case of the CAM noise process, we do not see noticeable changes in the density for $Y_j$ when we choose stepsize values smaller than
 $\delta t \approx \epsilon/10$ for simulating $y_t$.
 Therefore we simulate the fast CAM noise process with time discretization $\delta{t} \leq \epsilon/10$ throughout the paper. 

\section{Estimating $\Sigma$}
\label{sec:estim_sigma}
We estimate $\Sigma$ based on simulations of the integral of $ \int_{0}^{T}\,y_{s/\epsilon}\,ds$  with a value of $T= O(1)$ and a trapezoidal method, simulating $y_t$ in (\ref{eqn:CAMsde}) with a step size $\delta{t} \leq \epsilon/10$.  To ensure that the estimate for $\Sigma$ avoids the potential sensitivities described in Section \ref{discuss_num} we use a  value of  $\epsilon$ smaller than the values used in the numerical examples of Section  ($\epsilon = 10^{-5}$).  By choosing a small value of $\epsilon$, we seek an approximation for $\Sigma$ based on large $N_Y$.  As highlighted in Section \ref{discuss_num}, $N_Y$ increases with decreasing $\epsilon$ for $T$ fixed and $\Delta$  satisfying Conditions $A$ and $B$ in Section \ref{sec:reductionCAMtostable}.  
Then the distribution of $\int_{0}^{T}y_{s/\epsilon}\,ds$ is close to the $\alpha$-stable distribution (\ref{etaz_equiv}), and the approximation for $\Sigma$ is obtained from (\ref{etaz_equiv}) for given values of $\epsilon$ and $T$. 

The estimate of $\Sigma$ is based on a least squares fit of the  characteristic function of $ \int_{0}^{T}\,y_{s/\epsilon}\,ds = \sum_{j}Y_{j}$, $\psi_{S}(l;\sigma_Y) = \expected{\exp\left(il\int_{0}^{T}\,y_{s/\epsilon}\,ds\right)}$, to the empirical characteristic function 
\begin{equation}
\hat{\psi}_{S}(l) = \frac{1}{N_{S}}\sum_{j = 1}^{N_{S}}\exp\left(ilS_{j}^{(T)}\right), \quad S_{j}^{(T)} \sim \int_{0}^{T}\,y_{s/\epsilon}\,ds.
\end{equation}
based on realizations of $S_{j}^{(T)},\,j = 1,2,\dots,N_{S}$. The characteristic function $\psi_S$ is given in (\ref{eqn:apprx_charfuncY}),
\begin{equation}
\psi_S(l;\sigma_S) \to \exp\left(\sigma_S^{\alpha^{*}}|l|^{\alpha^{*}}\Xi(l,\alpha^{*},\beta^{*}) \right).
\end{equation}
The least squares fit of $\hat{\psi}_{S}(l)$ to $\psi_{S}(l;\sigma_S)$
is obtained by the minimization
\begin{equation}
\sigma_S= \operatornamewithlimits{\arg\min}_{\sigma}  \sum_{l = l_{\rm min}}^{l_{\rm max}}\left(\psi_{S}(l;\sigma) - \hat{\psi}_{S}(l) \right)^{2}.
\end{equation}
From the estimate of $\sigma_{S}$, we obtain our estimate for $\sigma_{Y} = \sigma_{S}/N_{Y}^{1/\alpha^{*}}$ where $\alpha^{*}$ is defined in (\ref{eqn:alpha_conv}). We note that the sum is proportional to integral estimates of the square of the difference between $\psi_{S}$ and $\hat{\psi}_{S}$ on the interval $[l_{\rm min}, l_{\rm max}]$. This numerical method is similar to the method of Koutrouvelis which estimates the parameters for $\alpha$-stable random variables via a least squares method and the characteristic function \cite{Koutrouvelis1980}. In our estimator, we use uniformly spaced points in the domain of $l$, but other studies have shown that a non-uniform spacing of points is optimal for such an estimator \cite{Besbeas2008}. While we opted for an unweighted estimator, there are also weighted versions of the Koutrouvelis method that may be worth considering. Nonetheless, our estimator is sufficient for obtaining accurate estimates of $\sigma_{Y}$ under the assumption that $S_{j}$ are distributed according to an $\alpha$-stable law. We use this estimator to verify the scaling relationships for $\Sigma$ with respect to $\Delta$ and $\epsilon$ (i.e. (\ref{sigma_Yprop}), as per Figure \ref{fig:sigma_dependence}), and estimate a value of $\theta$ via the relationship $\Sigma = \sigma^{*}/\theta$ (\ref{eqn:theta_true}) where $\sigma_{z} = \sigma^{*}$ to determine an equivalent OULp for the $y_{t}$ (as per Section \ref{sec:reductionCAMtostable}).

We use our estimator based on the empirical characteristic function $\hat\psi_{S}(l)$, rather than packages that can simultaneously estimate all parameters for a stable distribution such as the  {\tt STBL} package for MATLAB for values, which is based on the Koutrouvelis method \cite{STBL}. Tests with {\tt STBL} showed it to be unreliable in determining an estimate of the stability index $\alpha$ (which we know exactly) for values of $\alpha$ near 2, and consequently estimates of the corresponding scale parameter are not reliable.

\bibliography{biblio}

\end{document}